\documentclass{article}
\usepackage{tikz}
\usepackage{amsmath}
\usepackage{amsfonts}
\usepackage{amssymb}
\usepackage{xspace}
\usepackage{graphicx}
\date{}

\makeatletter \@addtoreset{equation}{section} \makeatother

\newcommand{\adl}{\vspace{1\baselineskip}}

\newenvironment{proof}{\par\noindent{\sc Proof:}
}{\hfill\llap{$\Box$}\vspace{1\baselineskip}\par\noindent}
\newenvironment{prooff}{\par\noindent{\sc Proof}
}{\hfill\llap{$\Box$}\vspace{1\baselineskip}\par\noindent}
\newtheorem{theorem}{Theorem}[section]
\newtheorem{proposition}[theorem]{Proposition}
\newtheorem{lemma}[theorem]{Lemma}
\newtheorem{corollary}[theorem]{Corollary}
\newtheorem{remark}[theorem]{Remark}
\newtheorem{definition}[theorem]{Definition}
\newtheorem{example}[theorem]{Example}
\newcommand{\beq}{\begin{equation}}
\newcommand{\eeq}{\end{equation}}

\newcommand{\ba}{\begin{array}}
\newcommand{\ea}{\end{array}}
\newcommand{\bt}{\begin{theorem}}
\newcommand{\et}{\end{theorem}}
\newcommand{\bp}{\begin{proposition}}
\newcommand{\ep}{\end{proposition}}
\newcommand{\bl}{\begin{lemma}}
\newcommand{\el}{\end{lemma}}
\newcommand{\bc}{\begin{corollary}}
\newcommand{\ec}{\end{corollary}}
\newcommand{\bi}{\begin{itemize}}
\newcommand{\ei}{\end{itemize}}
\newcommand{\ben}{\begin{enumerate}}
\newcommand{\een}{\end{enumerate}}
\newcommand{\bpf}{\begin{proof}}
\newcommand{\epf}{\end{proof}}
\newcommand{\bpff}{\begin{prooff}}
\newcommand{\epff}{\end{prooff}}
\newcommand{\bdf}{\begin{definition}\rm}
\newcommand{\edf}{\end{definition}}
\newcommand{\br}{\begin{remark}\rm}
\newcommand{\er}{\end{remark}}
\newcommand{\bex}{\begin{example}\rm}
\newcommand{\eex}{\end{example}}
\def\pri{\hbox to 10pt{\hfil\hbox to 0.4pt{\vrule height5pt width0.4pt
                 depth0pt}\vrule width5pt height0.4pt depth0pt\hfil}}
\newcommand{\TC}{{\rm TC}}

\newcommand{\TT}{{\rm TT}}

\newcommand{\gk}{{\bf{k}}}

\newcommand{\gK}{{\bf K}}

\newcommand{\gT}{{\bf T}}

\newcommand{\calL}{{\cal L}}

\newcommand{\gR}{{\mathbb R}}
\newcommand{\gZ}{{\mathbb Z}}

\newcommand{\Nat}{{\mathbb N}}

\newcommand{\RN}{{\mathbb R}^{N}}

\def\E{{\sf e}}
\def\IC{{\sf i}}

\newcommand{\om}{\omega}          

\newcommand{\D}{{\cal D}}
\newcommand{\F}{{\cal F}}

\newcommand{\calP}{{\cal P}}

\newcommand{\M}{{\cal M}}

\newcommand{\BV}{\mathop{\rm BV}\nolimits}
\newcommand{\SBV}{\mathop{\rm SBV}\nolimits}

\def\meas{\mathop{\rm meas}\nolimits}

\newcommand{\If}{\ \mbox{\rm if}\ }

\def\Lip{\mathop{\rm Lip}\nolimits}

\newcommand{\Sph}{{{\mathcal S}^2}}

\newcommand{\SP}{{{\mathbb S}^2}}
\newcommand{\Su}{{{\mathbb S}^1}}
\newcommand{\SN}{{{\mathbb S}^{N-1}}}

\newcommand{\gi}{\mathfrak{i}}
\def\mesh{\mathop{\rm mesh}\nolimits}
\newcommand{\ttt}{{\mathfrak T}}

\newcommand{\kk}{{\mathfrak K}}

\newcommand{\ee}{{\bf e}}
\newcommand{\gt}{{\bf t}}
\newcommand{\gn}{{\bf n}}
\newcommand{\gu}{{\bf u}}
\newcommand{\gv}{{\bf v}}

\newcommand{\gr}{{\bf r}}
\newcommand{\gc}{{\bf c}}

\newcommand{\Var}{\mathop{\rm Var}\nolimits}

\let\a=\alpha
\let\be=\beta

\let\e=\varepsilon
\let\f=\phi
\let\vf=\varphi
\let\g=\gamma

\let\m=\mu

\let\p=\pi

\let\vf=\varphi

\let\GG=\Gamma

\let\TT=\Theta

\let\wih=\widehat

\let\wid=\widetilde

\let\pa=\partial
\let\sb=\subset

\let\emp=\emptyset

\let\fa=\forall
\let\tim=\times

\let\sm=\setminus
\let\ol=\overline

\let\ds=\displaystyle

\let\lan=\langle
\let\ran=\rangle
\let\ii=\infty

\title{\Large \bf The total intrinsic curvature of curves in Riemannian surfaces \\ \& Erratum}
\author{\it Domenico Mucci and Alberto Saracco
\footnote{%
{\sc Dipartimento di Scienze Matematiche,
Fisiche ed Informatiche, Universit\`{a} di Parma,
Parco Area delle Scienze 53/A, I-43124 Parma, Italy.
E-mail: domenico.mucci@unipr.it, alberto.saracco@unipr.it}
}
}
\begin{document}
\topskip=1.5truecm \maketitle \topskip=1.5truecm \maketitle
\date{}
{\small {\bf Abstract.}
We deal with irregular curves contained in smooth, closed, and compact surfaces.
For curves with finite total intrinsic curvature, a weak notion of parallel transport of tangent vector fields is well-defined in the Sobolev setting.
Also, the angle of the parallel transport is a function with bounded variation, and its total variation is equal to an energy functional that depends on the ``tangential" component of the derivative of the tantrix of the curve.
We show that the total intrinsic curvature of irregular curves agrees with such an energy functional.
By exploiting isometric embeddings, the previous results are then extended to irregular curves contained in Riemannian surfaces.
Finally, the relationship with the notion of displacement of a smooth curve is analyzed.}
\adl\par\noindent
\adl\par\noindent
\adl\par\noindent

\section*{Erratum}
Our paper \cite{MSpal} appeared in \textit{Rendiconti del Circolo Matematico di Palermo}. After pubblication, while working on \cite{MSsph}, we realized that in the statements of the main results, Theorems from 1 to 9 and Propositon 3 (here Theorems from \ref{iTcomp} to \ref{Tequality} and Proposition~\ref{PGB}), one has to assume in addition that the curve $\gc$ is rectifiable. This Erratum will appear in  \textit{Rendiconti del Circolo Matematico di Palermo}. In this arXiv paper we added the appropriate hypothesis in said Theorems and Proposition.
\par The main point is that the equivalence in formula (2.7), namely:
%
$$ \TC_\M(\gc)<\ii\iff \TC(\gc)<\ii $$
holds true for rectifiable curves $\gc$, whereas it is false in general that if $\TC_\M(\gc)<\ii$, then $\TC(\gc)<\ii$.
If one e.g. takes a curve in $\Sph$, the unit sphere in $\gR^3$, that winds around an equator infinitely many times, its total intrinsic curvature is zero but its length and total curvature are both infinite.
\par Our mistake goes back to a flaw that we recently found in \cite[Thm. 6.3.1]{AR}, where Alexandrov-Reshetnyak erroneously stated that if the geodesic turn of a spherical curve is finite, then its spatial turn is also finite.
This is true if the spherical diameter of the curve is smaller than a dimensional constant $\delta_0$. In this case, in fact, for polygonal curves in $\Sph$ they obtain the inequality $\gk^*(P)\leq \p+2\gk_\Sph(P)$.

Therefore, their statement holds true provided that the curve can be divided in a finite number of arcs each one with spherical diameter smaller than $\delta_0$.
However, the latter property is false, in general, if the curve fails to be rectifiable, as the previous example shows.
\par Dealing with rectifiable curves $\gc$ in $\M$, in fact, by the smoothness and compactness of $\M$, the normal curvature of the geodesic arcs of $\M$ is uniformly bounded, and hence we recover the nontrivial implication $\Rightarrow$ in the previous equivalence by arguing as in the model case $\M=\Sph$ considered in \cite{AR}. 
\par For that reason, all the main results in \cite{MSpal} hold true for rectifiable curves with finite total intrinsic curvature.
\section{Introduction}
The theory of irregular curves goes back to A.~D.~Alexandrov and his collaborators in the 40's of the last century.
His joint work with Yu~G.~Reshetnyak is collected in the book \cite{AR} published in 1989.
We address to the survey paper \cite{Re} for detailed references.
\par A fundamental role in the theory of the Russian school is played by the class of {\em one-sidedly smooth} curves.
Such a regularity is exhibited e.g. by rectifiable curves in the Euclidean space $\RN$ with finite {\em total curvature}.
In fact, the unit tangent vector (or tantrix) exists almost everywhere, and it turns out to be a one-dimensional function of {\em bounded variation}.
%
By exploiting arguments based on integral geometric formulas, Alexandrov-Reshetnyak were also able to study irregular curves with values in the unit $N$-sphere.
\par A parallel theory of curves with finite total curvature, say {\sc ftc} curves, was introduced with a slightly different approach by J.~W.~Milnor \cite{Mi,Mi53} in the 50's.
More recently, J.~M.~Sullivan \cite{Su_curv} analyzed variational problems and geometric knot theory in this framework, showing the interplay between discrete and differential geometry.
For our purposes, we recall that the {\em total curvature} (i.e., the supremum of the {\em rotation} of the polygonals inscribed in the curve) of any {\sc ftc} curve in $\RN$ turns out to be equal to
the {essential variation} of the tantrix of the curve in the Gauss sphere $\SN$, see \eqref{VarSPt}.
For smooth curves, it clearly agrees with the integral of the scalar curvature.
\par Differently to the Euclidean case, an {\em intrinsic} theory of {\sc ftc} curves with values e.g. in a Riemannian manifold $\M$ fails to be complete,
even in the model case $\M=\Sph$, the unit 2-sphere in $\gR^3$.
\par A first problem comes with the good notion of {\em total intrinsic curvature} $\TC_\M(\gc)$ of an irregular curve $\gc$ in $\M$, in terms of
the best approximation with ``curved" polygonals of $\M$ inscribed in $\gc$.
In fact, for manifolds with positive sectional curvature (as e.g. $\M=\Sph$) the crucial monotonicity formula of the rotation of inscribed polygonals fails to hold.
\par In order to overcome this drawback, the good intrinsic notion turns out to be the one proposed by S.~B.~Alexander and R.~L.~Bishop \cite{AB}, that goes back to the one considered by Alexandrov-Reshetnyak \cite{AR}.
\par It involves the notion of {\em modulus} of an inscribed polygonal, that is, the greatest geodesic diameter of the arcs of the curve detected by the polygonal,
see Definition~\ref{Dcurv}.
\par With this notation, in fact, C.~Maneesawarng and Y.~Lenbury \cite{ML} showed that the total intrinsic curvature of a {\sc ftc} curve in $\M$ is equal to the limit of the rotation of {\em any} sequence of inscribed polygonals whose modulus goes to zero, see Proposition~\ref{Pappr}.
\par Notwithstanding, to our knowledge an explicit representation formula for the total intrinsic curvature $\TC_\M(\gc)$ is unknown in this general framework, for irregular curves $\gc$.
\par A partial result in this direction has been obtained by M.~Castrill\'on~Lopez, V.~Fernand\'ez~Mateos, and J.~Mu\~noz~Masqu\'e in \cite{CMM}
for the sub-class of (piecewise) smooth curves, see Theorem~\ref{Tdens}. Extending a result by Bishop \cite{Bi}, they showed that
\beq\label{ismooth} \TC_\M(\gc)=\int_\gc|\kk_g|\,ds +\sum_i|\a_i| \eeq
where $\kk_g$ is the {\em geodesic curvature} of the curve (that exists up to a finite number of points) and the second addendum is the finite sum of the ``turning angles" at the corner points of $\gc$.
\adl\par\noindent{\large\sc Content of the paper.} We deal with irregular curves contained in 2-dimensional Riemannian manifolds and with finite total intrinsic curvature. We first consider curves $\gc$ contained in a smooth (at least of class $C^3$), closed,
compact, and immersed surface $\M$ in $\RN$. Notice that $\M$ is not assumed to be oriented.
\par For the sake of clearness, in the first three sections we deal with the case of surfaces $\M$ in $\gR^3$, our model case being $\M=\Sph$, the standard unit sphere.
The high codimension case, $N\geq 4$, is treated in Sec.~\ref{Sec:N}.
\par We remark that the analysis of irregular curves in high dimension Riemannian manifolds needs some more work, and hence it will not be treated in this paper.
\par In Sec.~\ref{Sec:Tic}, we collect the notation concerning one dimensional $\BV$-functions, total curvature, geodesic curvature, and total intrinsic curvature,
by discussing the previously cited properties.
\par Our first new result, Theorem~\ref{Tcomp}, states that a notion of {\em weak parallel transport} is well-defined for curves with finite total intrinsic curvature.
For that reason, in Sec.~\ref{Sec:pt} we collect some well-known features concerning the classical parallel transport of tangent vector fields along smooth curves.
We also deal with piecewise-smooth curves, having in mind the case of polygonals $P_h$ in $\M$ inscribed in the irregular curve $\gc$.
\par Now, if the curve $\gc$ in $\M$ has finite total intrinsic curvature, say $\TC_\M(\gc)<\infty$, then $\gc$ is rectifiable.
We let $\gc:\ol I_L\to\M$ be its arc-length parameterization, where $I_L:=(0,L)$ and $L$ is the length of $\gc$.
By Rademacher's theorem, the tantrix $\gt:=\dot\gc$ is well-defined a.e. on $I_L$. Moreover,
by smoothness and compactness of $\M$, it turns out that $\gc$ is also a {\sc ftc} curve in $\RN$. Therefore, the tantrix $\gt$ is a function with bounded variation.
\par We also denote by $\gu$ the unit conormal to $\gc$ obtained by means of a positive rotation of $\gt$ on the tangent space $T_\gc\M$ along $\gc$.
If $\M\sb\gR^3$, we let $\gu:=\gn\tim\gt$, where $\gn$ is the (Lipschitz-continuous) outward unit normal to $\M$ along the curve.
\par In the sequel, the polygonals $P_h:\ol I_L\to\M$ are parameterized with constant velocity, and we denote by
$X_h:\ol I_L\to\RN$ the parallel transport of the vector field $\gt(0)$ along $P_h$.
Our Theorem~\ref{Tcomp} states:
\bt\label{iTcomp} If $\gc$ is a rectifiable curve, $\TC_\M(\gc)<\infty$, and
$\{P_h\}$ is a sequence of inscribed polygonals whose modulus goes to zero, then a subsequence of $\{X_h\}$ strongly converges in $W^{1,1}$ to some function
$X\in W^{1,1}(I_L,\RN)$ satisfying
$$ X(s)=\cos\Theta(s)\,\gt(s)-\sin\Theta(s)\,\gu(s) $$
for a.e. $s\in I_L$.
Furthermore, the {\em angle function} $\Theta$ has bounded variation, $\Theta\in\BV (I_L)$.
\et
\par For smooth curves $\gc$ on $\M$, the arc-length derivative $\dot\Theta$ of the angle function of the parallel transport is equal
to the geodesic curvature $\kk_g$ of the curve. In our second result, we shall compute the total variation of the
three components of the derivative of the optimal angle function $\Theta$, showing their relation with the three corresponding
components of the ``tangential derivative" of the tantrix $\gt:=\dot\gc$.
\par For this purpose, we recall that the distributional derivative of a $\BV$ function
$f:I_L\to\gR^k$ is a finite measure given by the sum $Df=D^af+D^Cf+D^Jf$ of its absolutely continuous, Cantor, and Jump components.
The latter ones are mutually singular and the decomposition $|Df|(I_L)=|D^af|(I_L)+|\D^Cf|(I_L)+|D^Jf|(I_L)$ of the total variation holds.
\par
The optimal angle is obtained by possibly minimizing the Jump of $\Theta$, without affecting the definition of weak parallel transport $X$, due to the $2\pi$-periodicity, see Remark~\ref{RTheta}.
Our Theorem~\ref{Tangle}, in fact, states:
\bt\label{iTangle} The optimal angle function $\Theta$ in Theorem $\ref{iTcomp}$ satisfies:
$$|D^a\Theta|(I_L)=\int_0^L|\dot\gt\bullet\gu|\,ds\,,\quad |D^C\Theta|(I_L)=|D^C\gt|(I_L)\,,\quad
|D^J\Theta|(I_L)=\sum_{s\in J_\gt}d_\SN(\gt(s+),\gt(s-)) $$
where $\bullet$ is the scalar product in $\RN$ and $\gt(s\pm)$ denotes the right or left limit of $\gt$ at $s$.
\et
\par As a consequence, the weak parallel transport $X$ along $\gc$ is essentially unique.
Notice, moreover, that for smooth curves $\gc$, in the first integral from Theorem~\ref{iTangle} one has $|\dot\gt\bullet\gu|=|\kk_g|$, whereas for piecewise smooth curves the Jump set $J_\gt$ of the tantrix is finite, and the last term (where $\gt(s\pm)$ denote the right and left limit of $\gt$ at the Jump points) agrees with the sum of the turning angles at the corner points.
\par For a curve $\gc$ with finite total intrinsic curvature, we are thus led to introduce the {\em energy functional}
\beq\label{iFt}
\F(\gt):=\int_0^L|\dot\gt\bullet\gu|\,ds+ |D^C\gt| (I_L)+\sum_{s\in J_\gt}d_\SN(\gt(s+),\gt(s-)) \eeq
where, we recall, $\gt:=\dot\gc$ is a function with bounded variation.
In the cited Theorem~\ref{Tdens} on piecewise smooth curves, in fact, formula \eqref{ismooth} reads:
\beq\label{iform} \TC_\M(\gc)=|D\Theta|(I_L)=\F(\gt)\,,\qquad \gt:=\dot\gc \eeq
\par We also point out that the Cantor component $D^C\gt$ of the derivative of the tantrix is tangential to $\M$.
More precisely, recalling that the unit conormal satisfies $\gu(s)\in T_{\gc(s)}\M$ for a.e. $s\in I_L$, we have:
$$ D^C\gt=\gu\,(\gu\bullet D^C\gt)=\gu\,D^C\Theta\,. $$
\par We thus expect that {\em the total intrinsic curvature $\TC_\M(\gc)$ agrees with the total variation $|D\Theta|(I_L)$ of the angle function},
and hence, by Theorem \ref{iTangle}, that the explicit formula \eqref{iform} holds true in full generality.
\par Now, denoting by $\Theta_h$ the angle function of the parallel transport $X_h$ along an approximating sequence $\{P_h\}$ as in Theorem \ref{iTcomp}, on
account of the cited Proposition~\ref{Pappr}, the representation formula \eqref{iform} holds true as a consequence of the {\em strict convergence}
\beq\label{istrict} \lim_{h\to\ii} |D\Theta_h|(I_L)=|D\Theta|(I_L)\,. \eeq
\par Obtaining the strict convergence \eqref{istrict} is a quite difficult task.
We observe that if one considers planar curves in $\gR^2$, the above limit holds true provided that one replaces the angle of the parallel transport with the oriented angle w.r.t. a fixed direction.
Therefore, in some sense, such a property relies on the validity of a ``planar" version of Gauss-Bonnet theorem, for domains whose boundary is parameterized by a curve with finite total curvature, see Sec.~\ref{Sec:GB}.
\par Following this approach, we show that the classical Gauss-Bonnet theorem generalizes to domains $U$ in $\M$ bounded by simple and closed
curves $\gc$ with finite total intrinsic curvature.
Referring to Theorem~\ref{TGB} for the precise statement, we only remark here that the term given by the circuitation of the geodesic curvature along the boundary of $U$, see \eqref{GB}, is replaced by the integral
$ \int_0^L k(s)\,ds$, where $k(s)\,ds:=D\Theta[0,s)$ and $\Theta$ is the angle function in Theorems~$\ref{iTcomp}$ and $\ref{iTangle}$, so that
$$ \int_0^L k(s)\,ds=\Theta(L)-\Theta(0)\,. $$
%
%
%
%
%
%
\par We point out that the class of curves with finite total intrinsic curvature seems to be the largest ambient in which the Gauss-Bonnet theorem makes sense.
If $\TC_\M(\gc)=\ii$, in fact, we expect that there is no way to find a finite measure that contains the information
(given by the derivative $D\Theta$ of the angle function of the parallel transport along the curve) on the ``signed geodesic curvature" of the curve $\gc$.
\par Our Lemma~\ref{LGB} on one-sidedly smooth curves, which is illustrated in Figure~$\ref{GBpic}$, allows to suitably exploit the generalized Gauss-Bonnet formulas from Theorem~\ref{TGB}.
In Proposition~\ref{PGB}, in fact, we build up a sequence $\{\wid\Theta_h\}$ of ``modified" angle functions that allows us to recover the upper semicontinuity
inequality in the strict convergence \eqref{istrict}, the lower semicontinuity inequality being a trivial matter.
We remark that a bit of care in the construction of the functions $\wid\Theta_h$ has to be taken when the surface $\M$ has positive Gauss curvature near the curve $\gc$, as in the model case $\M=\Sph$.
In conclusion, in Theorem~\ref{Tequality} we obtain:
\bt\label{iTequality} For every rectifiable curve $\gc$ in $\M$ with finite total curvature, $\TC_\M(\gc)<\ii$, the representation formula \eqref{iform} holds true,
where $\F(\gt)$ is the energy functional in \eqref{iFt} and $\gt=\dot\gc$ is the tantrix of the curve.
\et
\par In Sec.~\ref{Sec:Riem}, we deal with the case of curves into any smooth, closed, and compact Riemannian surface $\wid\M$.
The notion of total intrinsic curvature, in fact, clearly extends to curves $\g$ in $\wid\M$, where it is denoted by $\TC_{\wid\M}(\g)$.
\par By means of an isometric embedding $F$ of $\wid\M$ into a surface $\M=F(\wid\M)$ in $\RN$, we can apply our previous results to the curve $\gc:=F\circ\g$.
\par For this purpose, we shall focus in particular on the validity of the compactness theorem~\ref{iTcomp}.
In fact, by a quick inspection it turns out that the fundamental inequality \eqref{traspest} is the unique point of the previous theory where we used non-intrinsic quantities.
\par Moreover, we introduce geodesic polar coordinates, and write the local expression of the geodesic curvature of a smooth curve $\g$ in $\wid\M$.
It turns out that length, angles and geodesics are preserved by isometries.
In fact, we show that the geodesic curvature $\kk_g$ of $\gc:=F\circ\g$ in $\M:=F(\wid\M)$ agrees with the intrinsic local expression, and hence that the latter does not depend on the choice of the isometric embedding.
In a similar way, we check that the rotation of a polygonal $\wid P$ in $\wid\M$ is an intrinsic notion.
\par As a consequence, for piecewise smooth curves $\g$ in $\wid\M$ we obtain the equality:
$$ \TC_{\wid\M}(\g)=\TC_\M(\gc)\qquad \If\quad \gc:=F\circ\g $$
independently of the chosen isometric embedding $F$.
In conclusion, we obtain the following:
\bt\label{iTequalityS} For
every rectifiable curve $\g$ in $\wid\M$ with finite total intrinsic curvature, we have
$$ \TC_{\wid\M}(\g)=\F(\gt) $$
where the energy functional $\F(\gt)$ is defined by $\eqref{iFt}$ in correspondence to the tantrix $\gt=\dot\gc$ of $\gc=F\circ\g$, and $F$ is any isometric embedding of $\wid\M$ as above. \et
\par In Sec.~\ref{Sec:displ}, we finally deal with the notion of {\em development of a
smooth curve} $\g$ in a surface $\M$ of $\gR^3$, and analyze its relationship with the definition of total intrinsic curvature.
\par Namely, the {\em envelope of the tangent planes} to $\g$ is a ruled surface $\Sigma$
with zero Gauss curvature
around the trace of the curve, and hence it is locally isometric to a planar domain.
Moreover, the geodesic curvature $\kk_g$ of the curve $\g$ can be equivalently computed by using either local coordinates in $\M$ or in $\Sigma$.
\par The total intrinsic curvature $\TC_\Sigma(\g)$ of $\g$ as a curve in $\Sigma$ is well-defined, and
in Proposition \ref{Pdev} we show that it can
be recovered by means of the total curvature of the development of $\g$ in $\gR^2$, yielding to the expected formula:
$$\TC_\Sigma(\g)=\ds\int_\g|\kk_g|\,ds\,. $$
\par Therefore, even if in general the rotation of a polygonal $\wid P_h$ of $\Sigma$ and inscribed in $\g$, is different from the rotation of the corresponding polygonal $P_h$ in $\M$, see Example \ref{Edev}, by our previous results we infer that
$$ \TC_\M(\g)=\TC_\Sigma(\g) $$
which yields that the limits of the rotation of $P_h$ and of $\wid P_h$ coincide, if the modulus goes to zero.
\par We finally point out that similar arguments, based on considering iterations of the displacement of the ``complete tangent indicatrix", are proposed by Reshetnyak \cite{Re}
as a way to treat the ``curvatures" of an irregular curve in $\RN$.
A first step in this direction has been obtained in our paper \cite{MStor}, where a weak notion of torsion is analyzed.
\section{Total intrinsic curvature}\label{Sec:Tic}
In this section, we recall some properties concerning the total intrinsic curvature of smooth curves contained into surfaces.
We thus let $\M$ denote an immersed surface in $\gR^3$.
We assume $\M$ smooth (at least of class $C^3$), closed, and compact, our model case being $\M=\Sph$, the standard unit sphere in $\gR^3$.
%
\adl\par\noindent
{\large\sc $\BV$-functions of one variable.} We refer to Secs.~3.1 and 3.2 of \cite{AFP} for the following notation.
\par Let $I\sb\gR$ be a bounded open interval, and $N\in\Nat^+$.
A vector-valued summable function $u:I\to\RN$ is said to be of {\em bounded variation}
if its distributional derivative $Du$ is a finite $\RN$-valued measure in $I$.
%
\par
The {\em total variation} $|Du|(I)$ of a function $u\in\BV(I,\RN)$
is given by
$$ |Du|(I):=\sup\Bigl\{ \int_I \vf'(s)\,u(s)\,ds \mid \vf\in C^\infty_c(I,\RN)\,,\quad \Vert\vf\Vert_\infty\leq 1 \Bigr\} $$
and hence it does not depend on the choice of the representative in the
equivalence class of the functions that agree $\calL^1$-a.e. in $I$ with $u$, where
$\calL^1$ is the Lebesgue measure in $\gR$.
\par We say that a sequence $\{u_h\}\sb\BV(I,\RN)$ converges to $u\in \BV(I,\RN)$ {\em weakly-$^\ast$ in $\BV$} if $u_h$ converges to $u$ strongly in $L^1(I,\RN)$ and
$\sup_h|Du_h|(I)<\ii$. In this case, the lower semicontinuity inequality holds:
$$|Du|(I)\leq\liminf_{h\to\ii}|Du_h|(I)\,.$$
If in addition $|Du_h|(I)\to|Du|(I)$, we say that
$\{u_h\}$ {\em strictly converges} to $u$.
\par The {\em weak-$^\ast$ compactness} theorem yields that if $\{u_h\}\sb\BV(I,\RN)$ converges $\calL^1$-a.e. on $I$ to a function $u$,
and if $\sup_h|Du_h|(I)<\ii$, then $u\in\BV(I,\RN)$ and a subsequence of $\{u_h\}$ weakly-$^\ast$ converges to $u$.
%
\par Let $u\in\BV(I,\RN)$. Since each component of $u$ is the difference of two monotone functions, it turns out that $u$ is continuous outside an at most countable set, and that both the left and right limits
$u(s\pm):=\lim_{t\to s^\pm}u(t)$ exist for every $s\in I$. Also, $u$ is an $L^\ii$ function that is differentiable $\calL^1$-a.e. on $I$, with derivative $\dot u$ in $L^1(I,\RN)$.
\par The total variation of $u$ agrees with the {\em essential variation} $\Var_{\RN}(u)$, which is equal to the pointwise variation
of any {\em good representative} of $u$ in its equivalence class.
A good (or precise) representative is e.g. given by choosing $u(s)=(u(s+)+u(s-))/2$ at the discontinuity points.
Letting $u_\pm(s):=u(s\pm)$ for every $s\in I$, both the left- and right-continuous functions $u_\pm$ are good representatives.
\par If $u\in\BV(I,\RN)$, the decomposition into the {\em absolutely continuous}, {\em Jump}, and {\em Cantor} parts holds:
$$ Du=D^au+D^Ju+D^Cu\,,\quad |Du|(I)=|D^au|(I)+|D^Ju|(I)+|D^Cu|(I)\,.$$
More precisely, one splits $Du=D^au+D^su$ into the absolutely continuous and singular parts w.r.t. the Lebesgue measure $\calL^1$.
The Jump set $J_u$ being the (at most countable) set of discontinuity points of any good representative of $u$, and $\delta_s$ denoting the unit Dirac mass at $s\in I$,
one has:
$$D^a u=\dot u\,\calL^1\,,\qquad D^Ju=\sum_{s\in J_u}[u(s+)-u(s-)]\,\delta_s\,,\qquad D^Cu=D^su\pri(I\sm J_u)\,. $$
\par Also, any $u\in\BV(I,\RN)$ can be represented by $u=u^a+u^J+u^C$, where $u^a$ is a Sobolev function in $W^{1,1}(I,\RN)$, $u^J$ is a Jump function, and $u^C$ is a Cantor function, so that
$$ |D^au|(I)=|Du^a|(I)\,,\qquad |D^Ju|(I)=|Du^J|(I)\,,\qquad |D^Cu|(I)=|Du^C|(I)\,.  $$

\par Finally, we recall that if $u,v\in\BV(I):=\BV(I,\gR)$, the product $uv\in\BV(I)$.
In the particular case in which the Jump sets coincide, $J_u=J_v=J$, the chain rule formula (cf. \cite[Sec.~3.10]{AFP}) yields:
\beq\label{chain} D^a (uv)=(\dot uv+u\dot v)\,\calL^1\,,\quad D^J(uv)=\sum_{s\in J}[u(s+)v(s+)-u(s-)v(s-)]\,\delta_s\,,\quad D^C(uv)=u D^Cv+v D^Cu \eeq
where we can choose any good representatives of $u$ and $v$ in the third equality.
\adl\par\noindent
{\large\sc Total curvature.} We recall that the rotation $\gk^*(P)$ of a polygonal $P$ in $\gR^3$ is the sum of the exterior angles between consecutive segments. A polygonal $P$ is said to be inscribed in a curve $\gc:[a,b]\to\gR^3$ if $P$ is obtained by choosing a partition $a\leq t_0<t_1<\cdots < t_n\leq b$ and connecting with segments the consecutive points $\gc(t_i)$ of the curve. The mesh of the polygonal is $\mesh(P):=\max_{1\leq i\leq n}(t_i-t_{i-1})$. 
The {\em Euclidean total curvature} $\TC(\gc)$ of a curve $\gc$ in $\gR^3$ is defined by Milnor \cite{Mi,Mi53} as the supremum of the {\em rotation} $\gk^*(P)$ computed
among all the polygonals $P$ in $\gR^3$ which are inscribed in $\gc$.
Then $\TC(P)=\gk^*(P)$ for each polygonal $P$.
\par Let $\gc$ have compact support and finite total curvature, $\TC(\gc)<\ii$. Then, $\gc$ is a rectifiable curve.
In the sequel, we shall thus tacitly assume that $\gc$ is parameterized by arc-length, so that $\gc=\gc(s)$, with $s\in[0,L]=\ol I_L$, where
$I_L:=(0,L)$ and $L=\calL(\gc)$,
the length of $\gc$.
If $\gc$ is smooth and regular, one has $\TC(\gc)=\int_0^L|\gk|\,ds$,
where $\gk(s):=\ddot\gc(s)$ is the curvature vector.
More generally, since $\gc$ is a Lipschitz function, by Rademacher's theorem (cf. \cite[Thm.~2.14]{AFP}) it is differentiable $\calL^1$-a.e. in $I_L$.
Denoting by $\dot f:={d\over ds}f$ the derivative w.r.t.\ the arc-length parameter $s$, the tantrix $\gt=\dot\gc$ exists a.e., and actually
$\gt:I_L\to\gR^3$ is a function of bounded variation. Since moreover $\gt(s)\in\SP$ for a.e. $s$, where $\SP$ is the Gauss 2-sphere, we shall write
$\gt\in\BV(I_L,\SP)$.
The essential variation $\Var_\SP(\gt)$ of $\gt$ in $\SP$ differs from $\Var_{\gR^3}(\gt)$, as its definition involves the geodesic distance $d_\SP$ in $\SP$ instead of the Euclidean distance in $\gR^3$.
Therefore, $\Var_{\gR^3}(\gt)\leq \Var_\SP(\gt)$, and equality holds if and only if $\gt$ has a continuous representative.
More precisely, by decomposing $\gt=\gt^a+\gt^J+\gt^C$, one obtains:
\beq\label{VarSPt}  \Var_\SP(\gt)=\int_0^L|\dot\gt|\,ds+\sum_{s\in J_\gt}d_\SP(\gt(s+),\gt(s-))+|D^C\gt|(I_L) \eeq
whereas in the formula for $\Var_{\gR^3}(\gt)$, that is equal to $|D\gt|(I_L)$, one has to replace in \eqref{VarSPt} the geodesic distance $d_\SP(\gt(s+),\gt(s-))$ with the Euclidean distance
$|\gt(s+)-\gt(s-)|$ at each Jump point $s\in J_\gt$.
\adl\par\noindent
{\large\sc A representation formula.} The following facts hold:
\ben\item if $P$ and $P'$ are inscribed polygonals and $P'$ is obtained by adding a vertex in $\gc$ to the vertexes of $P$, then $\gk^*(P)\leq\gk^*(P')$\,;
\item if $\gc$ has finite total curvature, for each point $v$ in $\gc$, small open arcs of $\gc$ with an end point equal to $v$ have small total curvature. \een
\par As a consequence, compare \cite{Su_curv}, it turns out that
$\TC(\gc)=\Var_{\SP}(\gt)$, see \eqref{VarSPt}, and that the total curvature of $\gc$ is equal to the limit of $\gk^*(P_h)$ for {\em any} sequence $\{P_h\}$ of polygonals in $\gR^3$ inscribed in $\gc$ and such that $\mesh (P_h)\to 0$.
More precisely, if $\gt_h$ is the tantrix of $P_h$, then $\Var_\SP(\gt_h)\to\Var_\SP(\gt)$, see Remark~\ref{Rparall}.
%
%
\br\label{Rangle} The Cantor component $D^C\gt$ is non-trivial, in general. In fact, let e.g. $\g:\ol I\to\gR^2$, where $I=(0,1)$, denote the Cartesian curve $\g(t):=(t,u(t))$ in $\gR^2$
given by the graph of the primitive $u(t):=\int_0^tv(\lambda)\,d\lambda$
of the classical Cantor-Vitali function $v:\ol I\to \gR$
associated to the ``middle thirds" Cantor set. It turns out that $\gt=(1+v^2)^{-1/2}(1,v)$, whence $\gt$ is a Cantor function, i.e., $D^a\gt=D^J\gt=0$, and
$$D\gt(I)=D^C\gt(I)= \int_{I}{1\over(1+v^2)^{3/2}}\,(-v,1)\,d D^Cv\,. $$
Notice that the angle $\omega$ between the unit vectors $(1,0)$ and $\gt$ satisfies $\omega=\arctan v\in\BV(I)$. Therefore, $D\omega(I)=D^C\omega(I)=\int_{I}{1\over 1+v^2}\,d D^Cv $, which yields
$$ |D\omega|(I)=\int_{I}{1\over 1+v^2}\,d |D^Cv|=|D\gt|(I)=\TC(\g)={\p\over 4}\,. $$
\er
{\large\sc Geodesic curvature.}
Assume now that $\gc$ is a smooth and regular curve supported in $\M$.
The Darboux frame along $\gc$ is the triad $(\gt,\gn,\gu)$,
where $\gt(s):=\dot\gc(s)$ is the unit tangent vector, $\gn(s):=\nu(\gc(s))$, $\nu(p)$ being the unit normal to the tangent 2-space $T_p\M$, and $\gu(s):=\gn(s)\tim\gt(s)$, where $\tim$ denotes the vector product in $\gR^3$, is the unit conormal.
Therefore, the tangent space $T_{\gc(s)}\M$ is spanned by $(\gt(s),\gu(s))$.
%
The curvature vector $\gk(s)=\dot\gt(s)$ is orthogonal to
$\gt(s)$, and thus decomposes as
$$\gk(s)=\kk_g(s)\,\gu(s)+\kk_n(s)\,\gn(s) $$
where $\kk_g:=\gk\bullet\gu$ and $\kk_n:=\gk\bullet\gn$ denote the {\em geodesic} and {\em normal curvature} of $\gc$, respectively, and $\bullet$ is the scalar product in $\gR^3$.
The projection $\kk_g\gu$ of $\gk$ onto the tangent bundle of $\M$ is an intrinsic object, see Sec.~\ref{Sec:Riem}.
%
Also, the Frenet formulas in $\gR^3$ yield to the Darboux system:
\beq\label{Darboux} \dot\gt =\kk_g\gu+\kk_n\gn\,, \quad \dot\gn=-\kk_n\gt-\ttt_g\gu\,,\quad \dot\gu=-\kk_g\gt+\ttt_g\gn \eeq
where $\ttt_g:=\dot\gn\bullet(\gt\tim\gn)$ is the geodesic torsion of the curve.
%
\br\label{Rtors} If $\gc$ is a geodesic on $\M$, we have $\kk_g\equiv 0$, whence the Darboux frame $(\gt,\gn,\gu)$ agrees (up to the sign) with the Frenet frame, and the conormal $\gu$ with the bi-normal vector. In particular, the normal
curvature $\kk_n$ and the geodesic torsion $\ttt_g$ are equal (up to the sign) to the scalar curvature and to the torsion of $\gc$ in $\gR^3$, respectively.
Finally, the following estimate will be used in the proof of Theorem~\ref{Tcomp}: as for $\kk_n$, both $\ttt_g$ and its arc-length derivative are uniformly bounded by a constant only depending on the maximum of the modulus of the principal curvatures of $\M$ and of their derivatives, respectively.
\er
{\large\sc Total intrinsic curvature.} The {\em (intrinsic) rotation}
$\gk_\M(P)$ of a polygonal $P$ in $\M$, where $\M\sb\gR^3$, is the sum of the turning angles between the consecutive geodesic arcs of $P$.
The polygonal $P$ is  said to be inscribed in a curve $\gc:[a,b]\to\M\sb\gR^3$ if $P$ is obtained by choosing a partition $a\leq t_0<t_1<\cdots < t_n\leq b$ and connecting with geodesic segments the consecutive points $\gc(t_i)$ of the curve. 
For a general curve $\gc$ supported in $\M\sb\gR^3$, we shall denote by $\calP_\M(\gc)$ the class of polygonals in $\M$ which are inscribed in $\gc$. Also, if $\gc$ is rectifiable (and parameterized in arc-length) the mesh of a polygonal $P$ in $\calP_\M(\gc)$ is equivalently given by the maximum of the length of the arcs of $\gc$ bounded by the consecutive vertexes of $P$.
Notice that one clearly has $\gk_\M(P)\leq\TC(P)$, and that the difference $\TC(P)-\gk_\M(P)$ is equal to the sum of the integrals of the modulus of the
normal curvature $\kk_n$ of the geodesic arcs of $P$.
\par If e.g. $\M=\Sph$, then $\kk_n\equiv-1$ and hence $\TC(P)=\gk_\Sph(P)+\calL(P)$.
In general, by the smoothness and compactness of $\M$, the normal curvature of the geodesic arcs of $\M$ is uniformly bounded, and hence there exists a real constant $c_\M>0$
depending on $\M$ such that for each polygonal $P$ in $\M$
\beq\label{curvnorm} \TC(P)\leq\gk_\M(P)+c_\M\cdot\calL(P)\,. \eeq
\par
The following property has been proved in \cite{CMM}.
\bt\label{Tdens}{\bf (\cite[Thm.~3.4]{CMM})} Let $\gc$ be a regular curve in $\M$ of class $C^2$, parameterized by arc-length. Then, for any sequence $\{P_h\}\sb\calP_\M(\gc)$ such that $\mesh (P_h)\to 0$, one has
$$ \lim_{h\to\ii}\gk_\M(P_h)=\int_\gc|\kk_g|\,ds =\int_0^L|\kk_g(s)|\,ds\,. $$
\et
\par As a consequence, for a curve $\gc$ in $\M$, one is tempted to define its total intrinsic curvature as in the Euclidean case, i.e., as the supremum of the intrinsic rotation
$\gk_\M(P)$ computed among all the polygonals $P$ in $\calP_\M(\gc)$.
However, as observed in \cite{CMM}, if $\M$ has positive sectional (Gauss) curvature, as e.g. $\M=\Sph$, the latter definition does not work.
In fact, if $P,P'\in\calP_\M(\gc)$, and $P'$ is obtained by adding a vertex in $\gc$ to the vertexes of $P$, then the monotonicity inequality $\gk_\M(P)\leq\gk_\M(P')$ holds true in general provided that $\M$ has non-positive sectional curvature.
In fact, it relies on the fact that in this case the sum of the interior angles of a geodesic triangle of $\M$ is not greater than $\p$, see \cite[Lemma~4.1]{CMM}.
\bex\label{Emon}
If e.g. $\M=\Sph$, and $\gc$ is a parallel which is not a great circle, then the opposite inequality $\gk_\Sph(P)\geq\gk_\Sph(P')$ holds, and
for any $P\in\calP_\Sph(\gc)$ one has $\gk_\Sph(P)>\int_\gc|\kk_g|\,ds$, see Example~\ref{Eparal}. \eex
\par Actually, the good definition turns out to be the one introduced by Alexandrov-Reshetnyak \cite{AR}.
For this purpose, compare e.g. \cite{ML}, we recall that the {\em modulus} $\m_\gc(P)$ of a polygonal $P$ in $\calP_\M(\gc)$ is the maximum of the geodesic diameter of the arcs of $\gc$ determined by the consecutive vertexes in $P$.
For $\e>0$, we also let
$$ \Sigma_\e(\gc):=\{ P\in\calP_\M(\gc)\mid \m_\gc(P)<\e\}\,. $$
\bdf\label{Dcurv} The {\em total intrinsic curvature} of a curve $\gc$ in $\M$ is
$$ \TC_\M(\gc):=\lim_{\e\to 0^+}\sup\{ \gk_\M(P)\mid P\in \Sigma_\e(\gc)\}\,. $$ \edf
\par Clearly, the above limit is equal to the infimum of $\sup\{ \gk_\M(P)\mid P\in \Sigma_\e(\gc)\}$ as $\e>0$. Moreover, arguing as in \cite[Prop.~2.1]{ML}, for a polygonal $P$ in $\M$ we always have $\TC_\M(P)=\gk_\M(P)$.
Also, since $\M$ is compact, a curve with finite
total curvature $\TC_\M(\gc)<\ii$ is rectifiable, too (cf. \cite[Prop.~2.4]{ML}).
Most importantly, making use of a result by Dekster \cite{D}, as a consequence of \cite[Prop.~2.4]{ML} one obtains:
\bp\label{Pappr} The total curvature $\TC_\M(\gc)$ of any curve $\gc$ in $\M$ is equal to the limit of the rotation $\gk_\M(P_h)$ of {\em any} sequence of polygonals $\{P_h\}\sb\calP_\M(\gc)$ such that
$\m_\gc(P_h)\to 0$. \ep
\br\label{Rcurv} Proposition~\ref{Pappr} is proved in \cite[Thm.~6.3.2]{AR}, when $\M=\Sph$, and in \cite[Prop.~4.3]{CMM}, when $\M$ has non-positive Gauss curvature. The proof for general smooth surfaces $\M$ is obtained by arguing as in \cite[Prop.~2.4]{ML}, where it is firstly proved for curves in CAT(K) spaces. It suffices to observe that the Gauss curvature of $\M$ is bounded, provided that $\M$ is smooth and compact.
A crucial step is the following result (cf. \cite[Thm.~2.1.3]{AR}): if $\TC_\M(\gc)<\ii$, for each $\e>0$ there exists $\delta>0$ such that if $\gamma$ is an arc of $\gc$ with geodesic
diameter lower than $\delta$, the length of $\g$ is smaller than $\e$.
As a consequence, if $\{P_h\}\sb\calP_\M(\gc)$ is such that the modulus $\m_\gc(P_h)\to 0$, then also $\mesh(P_h)\to 0$, the converse implication being trivial. \er
\par Proposition~\ref{Pappr} fills the gap given by the lack of monotonicity observed e.g. in Example~\ref{Emon}, yielding to the conclusion that
Definition~\ref{Dcurv} involves a control on the modulus and not on the mesh, at least when the sectional curvature of $\M$ fails to be non-negative.
\par As a consequence, by Theorem~\ref{Tdens} one infers that for smooth curves $\gc$ in $\M$ one has $\TC_\M(\gc)=\int_\gc|\kk_g|\,ds$.
By \cite[Cor.~3.6]{CMM}, for piecewise smooth curves $\gc$ in $\M$ one similarly obtains that
\beq\label{TCMsmooth} \TC_\M(\gc)=\int_0^L|\kk_g(s)|\,ds+\sum_i|\a_i|\,. \eeq
In this formula, the integral is computed separately outside the corner points of $\gc$, where the geodesic curvature $\kk_g$ is well-defined,
and the second addendum denotes the finite sum of the absolute value of the oriented turning angles $\a_i$ between the incoming and outcoming unit tangent vectors at each corner point of $\gc$.
Therefore, for piecewise smooth curves we can rewrite formula \eqref{TCMsmooth} as
\beq\label{TCMsmoothV} \TC_\M(\gc)= \int_0^L|\dot\gt\bullet\gu|\,ds +
\sum_{s\in J_\gt}d_\SP(\gt(s+),\gt(s-))\,. \eeq
%
%
%
{\large\sc Properties.} For a curve $\gc$ in $\M$, we clearly have $\TC_\M(c)\leq\TC(\gc)$.
On account of the inequality \eqref{curvnorm}, arguing as in \cite[Thm.~6.3.1]{AR}, where the following property is proved for curves into $\Sph$, it turns out that if $\TC_\M(\gc)<\ii$, then also $\TC(\gc)<\ii$, and hence that we definitely have:
\beq\label{TCTCM}
 \TC_\M(\gc)<\ii\iff \TC(\gc)<\ii\,. \eeq
\par Therefore, if $\TC_\M(\gc)<\ii$, then $\gc$ is rectifiable and the tantrix $\gt:=\dot\gc\in\BV(I_L,\SP)$.
Moreover, the curve is {\em one-sidedly smooth} in the sense of \cite[Sec.~3.1]{AR}, i.e., the curve has a left and a right tangent $\gT_\pm(s)$ at all the points
$\gc(s)$ in the ``strong sense".
\br\label{Rossm} This implies that for each $s\in[0,L[$ and $\delta>0$ we can find $\e>0$ such that any secant inscribed in the arc $\gc_{\vert\,[s,s+\e]}$ forms with the straight line
$\gT_+(s)$ an angle less than $\delta$, and similarly for the left tangent. \er
\par
As in the smooth case, we let $\gn:=\nu\circ \gc$ denote the unit normal to $T_\gc\M$ along $\gc$. Since $\M$ is smooth and compact, and $\gc$ is Lipschitz-continuous, it turns out that
$\gn\in \Lip( [0,L],\SP)$.
Therefore, the {\em weak conormal} $\gu:=\gn\tim\gt$ belongs to $\BV(I_L,\SP)$, with $J_\gu=J_\gt$.
Since moreover $\dot\gt\bullet \gt=0$ a.e. in $I_L$, we may decompose
$\dot\gt=(\dot\gt\bullet\gu)\,\gu+(\dot\gt\bullet\gn)\,\gn$.
\br\label{RCantor} We finally see that if $\gc$ is a curve in $\M$ with finite total curvature, the Cantor component $D^C\gt$ is tangential to $\M$, namely:
$$ D^C\gt=\gu\,(\gu\bullet D^C\gt) $$
where $\gu(s)\in T_{\gc(s)}\M$ for a.e. $s\in I_L$.
In fact, using that $|\gt|^2=|\gu|^2=1$ and $\gt\bullet\gu=\gt\bullet\gn=0$ a.e., whereas both $\gt\bullet\gu$ and $\gt\bullet\gn$ are functions of bounded variation, and $D^C\gn=0$, by \eqref{chain} we infer that $\gt\bullet D^C\gt=0$, $\gu\bullet D^C\gu=0$,
$\gu\bullet D^C\gt=-\gt\bullet D^C\gu$, and $\gn\bullet D^C\gt=D^C(\gt\bullet\gn)=0$. Since $(\gt,\gn,\gu)$ is an orthonormal frame to $\gR^3$, the tangential property follows.
\er
\section{Parallel transport}\label{Sec:pt}
In this section, we collect some well-known facts concerning the parallel transport of tangent vector fields $X$ along smooth curves in $\M$. We then also analyze the case of piecewise smooth curves.
Finally, we give some more detail in the model case $\M=\Sph$.
\par Let $\gc$ be a smooth, regular, and rectifiable curve in $\M$. Then $X: [0,L]\to \gR^3$ is a parallel transport along $\gc$ if for each $s\in [0,L]$ one has $X(s)\in T_{\gc(s)}\M$ and
$\dot X(s)\perp T_{\gc(s)}\M$, i.e., $\dot X(s) \parallel\gn(s)$.
We recall that since ${d\over ds}|X(s)|^2=2 X(s)\bullet \dot X(s)=0$ for every $s$, the parallel transport preserves the length of the initial tangent vector $X(0)$.
\par The proof of the following well-known property is taken from \cite[13.6.1]{Pre}.
\bp\label{PTheta} Let $\Theta(s)$ denote the oriented angle from the parallel transport $X(s)$ to the tangent vector $\gt(s)$ to $\gc$. Then, the geodesic curvature of $\gc$ satisfies $\kk_g(s)=\dot\Theta(s)$ for each $s\in [0,L]$.
\ep
\bpf Assume $|X(0)|=1$, so that $|X(s)|=1$ for every $s$. Writing
\beq\label{Theta} X(s)=\cos\Theta(s)\,\gt(s)-\sin\Theta(s)\,\gu(s)\,,\qquad s\in [0,L] \eeq
we find for each $s$
$$ 0=\gt\bullet\dot X=\gt\bullet [(\cos\Theta\,\dot\gt-\sin\Theta\,\dot\gu)-\dot\Theta\,(\sin\TT\,\gt+\cos\TT\,\gu) ]=-\sin\TT\,(\gt\bullet\dot\gu+\dot\TT) $$
where we used that $\gt\bullet\dot\gt=\gt\bullet\gu=0$. Similarly, condition $\gu\bullet\dot X=0$ implies
$$ 0=\cos\TT\,(\gu\bullet\dot\gt-\dot\TT)\,. $$
Since $\gk=\dot\gt$, we have $\kk_g=\dot\gt\bullet\gu$. Using that $\gt\bullet \gu=0$, we also get
$\gt\bullet\dot\gu=-\dot\gt\bullet\gu=-\kk_g\,. $
Therefore, the above centered equations become
$$ (\kk_g(s)-\dot\TT(s))\,\sin\TT(s)=0=(\kk_g(s)-\dot\TT(s))\,\cos\TT(s)\qquad \fa\, s\in[0,L]$$
which yields $\kk_g=\dot\TT$. \epf
\par We thus get the formula for the total intrinsic curvature of a smooth regular curve $\gc$ in $\M$
\beq\label{TCM}
 \TC_{\M}(\gc)=\int_0^L|\kk_g(s)|\,ds=\int_0^L|\dot\TT(s)|\,ds  \eeq
compare e.g. \cite{CMM}. Finally, notice that when $X(s)\bullet\gt(s)\neq 0$, by \eqref{Theta} one has
\beq\label{Theta2} \tan\Theta(s)=-\frac{X(s)\bullet\gu(s)}{X(s)\bullet\gt(s)}\,.  \eeq
{\large\sc Piecewise smooth curves.}
The parallel transport \eqref{Theta} is a well-defined smooth vector field for each regular and piecewise smooth curve $\gc$, once the initial position $X(0)$ is prescribed. If e.g. the curve is rectifiable and its arc-length parameterization is piecewise $C^k$, then the parallel transport is of class $C^k$.  
Moreover, the angle $\Theta$ is a function of bounded variation, with a finite number of Jump points in correspondence to the values $\{s_i\mid i=1,\ldots n\}$
of the arc-length parameter $s\in I_L$ where $\gc(s)$ fails to be smooth, the corner points $\gc(s_i)$ of $\gc$.
More precisely, $\Theta$ is a special function of bounded variation in $\SBV(I_L)$, i.e., $D^C\Theta=0$, and its distributional derivative decomposes as
$ D\Theta =\dot\Theta\,\calL^1+D^J\Theta\,. $
\par By Proposition~\ref{PTheta}, it turns out that the derivative $\dot\Theta$ agrees with the geodesic curvature $\kk_g$ outside the corner points of $\gc$, and the Jump component $D^J\Theta$ is a sum of Dirac masses centered at the points $s_i$, with weight given by the oriented turning angles $\a_i$
between the incoming and outcoming unit tangent vectors at each corner point of $\gc$. We thus have
$$ D\Theta = \kk_g\,\calL^1+\sum_{i=1}^n\a_i\,\delta_{s_i}\,,\qquad |D\Theta| (I_L)=\int_0^L|\kk_g|\,ds+\sum_{i=1}^n|\a_i| $$
and hence by \eqref{TCMsmooth} one infers that
$$ |D\Theta| (I_L)=\TC_\M(\gc)\,. $$
In particular, if $\gc$ is a polygonal $P$ in $\M$, the angle function is piecewise constant and
$$ D\Theta = \sum_{i=1}^n\a_i\,\delta_{s_i}\,,\qquad |D\Theta| (I_L)=\sum_{i=1}^n|\a_i|=\gk_\M(P)\,.  $$

\par Moreover, denoting by $(\gt,\gn,\gu)$ the Darboux frame of $\gc$, so that formulas \eqref{Darboux} hold true outside the points $s_i$,
by the smoothness of $X$ in general we have
$$ \dot X=
-\sin\Theta\,\dot\Theta\,\gt-\cos\Theta\,\dot\Theta\,\gu+ \cos\Theta\,\dot\gt-\sin\Theta\,\dot\gu $$
and hence the parallel transport of piecewise smooth curves satisfies, for $s\neq s_i$,
\beq\label{dotX} \dot X= (\cos\Theta\,\kk_n-\sin\Theta\,\ttt_g)\,\gn\,. \eeq
{\large\sc Curves into the 2-sphere.} Assume now $\M=\Sph$.
Taking polar coordinates
$$ \gr(\theta,\vf)^T=(\sin\theta\cos\vf,\sin\theta\sin\vf,\cos\theta)\,,\qquad  \theta\in[0,\p]\,,\quad \vf\in[0,2\p] $$
the curve $\gc$ may thus be parameterized by $\gc(s)=\gr(\theta(s),\vf(s))^T$ for some smooth angle functions $\theta(s)$ and $\vf(s)$.
Consider the frame
$$ \ee_\theta(\theta,\vf):= \left(\begin{array}{c}\cos\theta\cos\vf \\ \cos\theta\sin\vf \\ -\sin\theta \\  \end{array}\right)\,,\qquad
\ee_\vf(\theta,\vf):= \left(\begin{array}{c}-\sin\vf \\ \cos\vf \\ 0\\  \end{array}\right)\,,\qquad
\gn(\theta,\vf):= \left(\begin{array}{c}\sin\theta\cos\vf \\ \sin\theta\sin\vf \\ \cos\theta \\  \end{array}\right)
$$
where $\gn=\ee_\theta\tim\ee_\vf$ is the outward unit normal. The partial derivatives of the tangent frame $(\ee_\theta,\ee_\vf)$ satisfy
$$ \pa_\theta\ee_\theta=-\gn\,,\quad \pa_\vf\ee_\theta=\cos\theta\,\ee_\vf\,,\quad \pa_\theta\ee_\vf\equiv 0\,,\quad \pa_\vf\ee_\vf=-\sin\theta\,\gn-\cos\theta\,\ee_\theta\,.$$
Letting
$$\ee_\theta(s):=\ee_\theta(\theta(s),\vf(s))\,,\quad \ee_\vf(s):=\ee_\vf(\theta(s),\vf(s))\,, \quad \gn(s):=\gn(\theta(s),\vf(s)) $$
we thus have
\beq\label{tang} \gt(s):=\dot\gc(s)=\dot\theta(s)\,\ee_\theta(s)+\sin\theta(s)\,\dot\vf(s)\,\ee_\vf(s)\,,\qquad\dot\theta(s)^2+\sin^2\theta(s)\,\dot\vf(s)^2=1\qquad\fa\,s\in[0,L]\,. \eeq
\par Consider a tangent vector field $X$ along $\gc$, so that
$$ X(s):=\a(s)\,\ee_\theta(s)+\be(s)\,\ee_\vf(s)\,,\qquad s\in[0,L] $$
for some smooth unknown functions $\a(s)$ and $\be(s)$.
We compute for each $s\in[0,L]$
$$ \ba{rl}\dot X= &\dot \a\,\ee_\theta+\a\,(\pa_\theta\ee_\theta\,\dot\theta+\pa_\vf\ee_\theta\,\dot\vf)+\dot \be\,\ee_\vf+\be\,(\pa_\theta\ee_\vf\,\dot\theta+\pa_\vf\ee_\vf\,\dot\vf) \\
= & \dot \a\,\ee_\theta+\a\,(-\gn\,\dot\theta+\cos\theta\,\ee_\vf\,\dot\vf)+\dot \be\,\ee_\vf+\be\,(-\sin\theta\,\gn\,\dot\vf-\cos\theta\,\ee_\theta\,\dot\vf) \\
= & (\dot\a-\be\,\cos\theta\,\dot\vf)\,\ee_\theta+(\dot\be+\a\,\cos\theta\,\dot\vf)\,\ee_\vf+(-\a\,\dot\theta-\be\,\sin\theta\,\dot\vf)\,\gn\,.
\ea$$
Condition for a parallel transport is $\dot X(s) \parallel\gn(s)$ for each $s$. This is equivalent to the first order system for the unknown coefficients $\a(s)$ and $\be(s)$\,:
\beq\label{ode1} \left\{ \ba{ll} \dot\a(s)=\cos\theta(s)\,\dot\vf(s)\,\be(s) \\ \dot\be(s)=-\cos\theta(s)\,\dot\vf(s)\,\a(s) \ea \right. \qquad s\in[0,L]
\eeq
which turns out to have a unique solution for any given initial position $X(0)\in T_{\gc(0)}\Sph$.
\par Since the parallel transport preserves the length, assuming $X(0)=\gt(0)$, we have
$$\a^2(s)+\be^2(s)=1\qquad\fa\,s\in[0,L]\,. $$
Therefore, from \eqref{ode1} one also obtains the identity:
\beq\label{idtrans}
\dot\a(s)\,\be(s)-\a(s)\,\dot\be(s)=\cos\theta(s)\,\dot\vf(s)\qquad \fa\,s\in[0,L]\,. \eeq
\par On account of \eqref{Theta2}, and since by \eqref{tang} the unit conormal along $\gc$ is
\beq\label{uu}
 \gu(s):=\gn(s)\tim\gt(s)= -\sin\theta(s)\,\dot\vf(s)\,\ee_\theta(s)+\dot\theta(s)\,\ee_\vf(s) \eeq
one infers that for each $s\in[0,L]$ such that $\a\,\dot\theta+\be\,\sin\theta\,\dot\vf\neq 0$,
$$ \tan \Theta=\frac{\a\,\sin\theta\,\dot\vf-\be\,\dot\theta}{\a\,\dot\theta+\be\,\sin\theta\,\dot\vf}\,. $$
Using repeatedly that $\a^2+\be^2\equiv \dot\theta^2+\sin^2\theta\,\dot\vf^2\equiv 1$, one has
$$ \ba{rl} \dot\Theta= & \ds{d\over ds}\,\bigl(\a\,\sin\theta\,\dot\vf-\be\,\dot\theta \bigr)\cdot\bigl( \a\,\dot\theta+\be\,\sin\theta\,\dot\vf \bigr)-
{d\over ds}\,\bigl( \a\,\dot\theta+\be\,\sin\theta\,\dot\vf \bigr)\cdot\bigl(\a\,\sin\theta\,\dot\vf-\be\,\dot\theta \bigr) \\
= & \ds \dot\a\,\be-\a\,\dot\be+\sin\theta\,(\ddot\vf\,\dot\theta-\ddot\theta\,\dot\vf)+\cos\theta\,\dot\theta^2\,\dot\vf \\
= & \sin\theta\,(\ddot\vf\,\dot\theta-\ddot\theta\,\dot\vf)+\cos\theta\,\dot\vf\,(\sin^2\theta\,\dot\vf^2+2\dot\theta^2)  \ea$$
where the last equality follows from the identity \eqref{idtrans}.
\par On the other hand, recalling formula \eqref{tang}, the curvature vector of $\gc$ is
\beq\label{curvS} \gk=\dot\gt=(\ddot\theta-\sin\theta\,\cos\theta\,\dot\vf^2)\,\ee_\theta+(2\cos\theta\,\dot\theta\,\dot\vf+\sin\theta\,\ddot\vf)\,\ee_\vf-\gn \eeq
and hence by \eqref{uu} the geodesic curvature becomes
\beq\label{geodcurvS} \kk_g=\gk\bullet\gu =\sin\theta\,(\ddot\vf\,\dot\theta-\ddot\theta\,\dot\vf)+\cos\theta\,\dot\vf\,(\sin^2\theta\,\dot\vf^2+2\dot\theta^2) \eeq
where $(\sin^2\theta\,\dot\vf^2+2\dot\theta^2)=(1+\dot\theta^2)$, so that one recovers the equality $\kk_g=\dot\Theta$ from Proposition~\ref{PTheta}.
\bex\label{Eparal} If $\gc=\gc_{\theta_0}$ is the parallel with constant co-latitude $\theta_0\in ]0,\p/2]$, we choose
$\theta(s)\equiv\theta_0$ and $\vf(s)=s/\sin\theta_0$, where $s\in[0,L]$, with $L:=\calL(\gc_{\theta_0})=2\p\,\sin\theta_0$. By
\eqref{tang} and \eqref{uu}, one has
$$ \gt(s)=\ee_\vf(\theta_0,s/\sin\theta_0)\,,\quad \gu(s)=-\ee_\theta(\theta_0,s/\sin\theta_0)
\qquad \fa\, s $$
and by solving the system \eqref{ode1} as above, on account of \eqref{curvS} and
\eqref{geodcurvS} one obtains
$$ \Theta(s)=\cot\theta_0\cdot s\,,\qquad \kk_g=\dot\Theta\equiv \cot\theta_0\qquad \fa\,s\,. $$
Therefore, according to \eqref{TCM} one recovers for any $\theta_0\in]0,\p/2]$ the formula
$$ \TC_{\Sph}(\gc_{\theta_0})=\int_0^{2\p\,\sin\theta_0}|\dot\Theta(s)|\,ds=2\p\,\cos\theta_0 $$
for the total intrinsic curvature of the parallel, compare e.g. \cite{CMM}. In particular, $\TC_\Sph(\gc_{\theta_0})$ is equal to zero when $\theta_0=\p/2$, i.e., when $\gc_{\theta_0}$ is a
great circle, whence a geodesic in $\Sph$.
\eex
\section{Weak parallel transport}
In this section, we show that a weak notion of parallel transport holds true for curves $\gc$ in $\M$ with finite total intrinsic curvature, see Theorem~\ref{Tcomp}.
The parallel transport turns out to be a Sobolev function satisfying \eqref{Theta}, where the unit tangent $\gt$ and conormal $\gu$ are functions of bounded variation, and the angle function $\Theta$ is of bounded variation, too.
As a consequence, we infer that {\em the optimal angle function $\Theta$ is essentially unique}, and that the {\em weak transport} $X$ along the non-smooth curve $\gc$
is well-defined by the $W^{1,1}$ tangent vector field in Theorem~\ref{Tcomp}.
In fact, it turns out that the distributional derivative of the angle function $\Theta$ is strongly related to the tangential component of the derivative of the tantrix $\gt$,
see Theorem~\ref{Tangle}.
\adl\par\noindent
{\large\sc A compactness result.} We first prove the following
\bt\label{Tcomp} Let $\gc$ be a rectifiable curve in $\M$ with finite total intrinsic curvature, parameterized by arc-length $\gc:[0,L]\to\M$, with $L=\calL(\gc)$. Let $\{P_h\}\sb\calP_\M(\gc)$ be such that the modulus $\m_\gc(P_h)\to 0$. For each $h$, let $P_h:[0,L]\to\M$ be parameterized with constant velocity, and let
$X_h:[0,L]\to\gR^3$ be the parallel transport along $P_h$, with constant initial condition $X_h(0)=\gt(0)\in\SP$. Then, possibly passing to a subsequence, the sequence $\{X_h\}$ strongly converges in $W^{1,1}$ to some function
$X\in W^{1,1}( I_L,\gR^3)$ satisfying
\beq\label{weakX} X(s)=\cos\Theta(s)\,\gt(s)-\sin\Theta(s)\,\gu(s) \eeq
%
for $\calL^1$-a.e. $s\in I_L$, where $\gt=\dot \gc$ is the unit tangent vector, $\gn$ the normal to $T_\gc\M$ along $\gc$, and $\gu:=\gn\tim\gt$ is the unit conormal.
Furthermore, $\gt$ and $\gu$ are functions in $\BV( I_L,\SP)$, and the angle function $\Theta$ has bounded variation in $\BV (I_L)$.
\et
\bpf Write for each $h$
\beq\label{Xh}  X_h(s)=\cos\Theta_h(s)\,\gt_h(s)-\sin\Theta_h(s)\,\gu_h(s) \eeq
and recall that $|D\Theta_h| (I_L)=\gk_\M(P_h)$, whereas the difference $\TC(P_h)-\gk_\M(P_h)$ is equal to the sum of the integrals of the modulus of the
normal curvature $\kk_n$ of the geodesic arcs of $P_h$, so that the inequality \eqref{curvnorm} holds. Using that $\gk_\M(P_h)\to \TC_\M(\gc)<\ii$ and
$\calL(P_h)\to L$, we thus obtain the bounds:
$$ \sup_h|D\Theta_h| (I_L)<\ii\,,\qquad \sup_h\Var_\SP(\gt_h)=\sup_h\TC(P_h)<\ii\,. $$
Therefore, by the weak-$^*$ compactness, and by using the strong convergence of $P_h$ to $\gc$, possibly passing to a subsequence it turns out that $\{\gt_h\}$ and $\{\gu_h\}$ converge weakly-$^\ast$ in the $\BV$-sense to $\gt$ and $\gu$, respectively, and that the sequence $\{\Theta_h\}$ converges weakly-$^\ast$ in the $\BV$-sense to
some function $\Theta\in\BV (I_L)$.
\par We claim that for each $s\in[0,L]$ and for $\delta>0$ small
\beq\label{traspest}
\int_0^L |\dot X_h(s+\delta)-\dot X_h(s)|\,ds  \leq C_\M\cdot\delta\cdot[\calL(P_h)+|D\Theta_h| (I_L) ] \eeq
where the real constant $C_\M$ only depends on $\M$.
As a consequence, the sequences $\{\calL(P_h)\}$ and $\{|D\Theta_h| (I_L)\}$ being bounded, it turns out that
$$ \lim_{|\delta|\to 0}\sup_h\int_0^L |\dot X_h(s+\delta)-\dot X_h(s)|\,ds=0 $$
whereas $|X_h(s)|\equiv 1$ for each $h$. Therefore, by Kolmogorov-Riesz-Frech\'et compactness theorem, a further subsequence of $\{X_h\}$ strongly converges in $W^{1,1}$ to some function
$X\in W^{1,1}( I_L,\gR^3)$. Finally, by the $L^1$ convergence of $\gt_h$, $\gu_h$ and $\Theta_h$ to $\gt$, $\gu$, and $\Theta$, respectively, we conclude that
\eqref{weakX} holds $\calL^1$-a.e. on $I_L$.
\par In order to prove the inequality \eqref{traspest}, for each $h$ we first smoothly extend the transport $X_h$ to an interval $[-\delta_0,L+\delta_0]$ along the extreme geodesic arcs of $P_h$, where $\delta_0>0$ is fixed. For $0<|\delta|<\delta_0$, using formula \eqref{dotX} for $X=X_h$,
and omitting for simplicity to write the dependence on $h$, for each $s\in[0,L]$ we have:
$$ \ba{rl} \dot X(s+\delta)-\dot X(s)= & (\cos\Theta(s+\delta)-\cos\Theta(s))\,\kk_n(s+\delta)\,\gn(s+\delta) \\
& +\cos\Theta(s)\,(\kk_n(s+\delta)-\kk_n(s))\,\gn(s+\delta) \\
& + \cos\Theta(s)\,\kk_n(s)\,(\gn(s+\delta)-\gn(s)) \\
& - (\sin\Theta(s+\delta)-\sin\Theta(s))\,\ttt_g(s+\delta)\,\gn(s+\delta) \\
& -\sin\Theta(s)\,(\ttt_g(s+\delta)-\ttt_g(s))\,\gn(s+\delta) \\
& - \sin\Theta(s)\,\ttt_g(s)\,(\gn(s+\delta)-\gn(s))\,.
\ea $$
\par On account of Remark~\ref{Rtors}, we first estimate the three terms depending on $\kk_n$ as follows:
$$ |(\cos\Theta(s+\delta)-\cos\Theta(s))\,\kk_n(s+\delta)\,\gn(s+\delta)|\leq |\kk_n(s+\delta)|\cdot |D\Theta|(s,s+\delta)\,, $$
where by Fubini-Tonelli's theorem
$$ \int_0^L|\kk_n(s+\delta)|\cdot |D\Theta|(s,s+\delta)\,ds\leq c_\M\cdot|D\Theta| (I_L)\cdot\delta\,, $$
$c_\M$ being the maximum of the modulus of the principal curvatures of $\M$.
Moreover,
$$ |\cos\Theta(s)\,(\kk_n(s+\delta)-\kk_n(s))\,\gn(s+\delta)|\leq c'_\M\,\delta $$
$c'_\M$ being the maximum of the modulus of the derivative of the principal curvatures of $\M$, and similarly, since
$|\gn(s+\delta)-\gn(s)|\leq {c_\M}\cdot \delta$, we get:
$$ |\cos\Theta(s)\,\kk_n(s)\,(\gn(s+\delta)-\gn(s))|\leq {c_\M}^2 \delta\,. $$
\par As to the three terms depending on $\ttt_g$, we infer as above:
$$ \int_0^L|(\sin\Theta(s+\delta)-\sin\Theta(s))\,\ttt_g(s+\delta)\,\gn(s+\delta)|\,ds\leq K_\M\cdot|D\Theta| (I_L)\cdot\delta $$
$K_\M$ being a uniform bound, only depending on $\M$, of the maximum of the modulus of the geodesic torsion of $P_h$, outside the corner points.
Moreover,
$$ |\sin\Theta(s)\,(\ttt_g(s+\delta)-\ttt_g(s))\,\gn(s+\delta)|\leq K'_\M\,\delta $$
$K'_\M$ being a uniform bound, only depending on $\M$, of the maximum of the modulus of the derivative of the geodesic torsion of $P_h$, outside the corner points. Finally,
$$ |\sin\Theta(s)\,\ttt_g(s)\,(\gn(s+\delta)-\gn(s))|\leq {K_\M}\,c_\M\, \delta\,. $$
Therefore, inequality \eqref{traspest} readily follows, and the proof is complete. \epf
{\large\sc The angle function.} In principle, the angle function $\Theta$ depends on the subsequence corresponding to the approximating sequence $\{P_h\}$.
We now show that the {\em optimal} angle function $\Theta$, see Remark~\ref{RTheta}, is essentially unique and hence that the parallel transport $X$ along irregular curves
$\gc$ with finite total curvature is well-defined in the $W^{1,1}$ setting. In fact, in Theorem~\ref{Tangle} we write the total variation of the optimal angle function in terms of the tangential weak derivative of the tantrix $\gt$.
\par
For this purpose, recalling the decomposition $\dot\gt=(\dot\gt\bullet\gu)\,\gu+(\dot\gt\bullet\gn)\,\gn$ of the differential of the tantrix $\gt:=\dot\gc$ into the tangential and normal components, we introduce the {\em energy functional}
\beq\label{Ft} \F(\gt):=\int_0^L|\dot\gt\bullet\gu|\,ds+ |D^C\gt| (I_L)+\sum_{s\in J_\gt}d_\SP(\gt(s+),\gt(s-))\,. \eeq
\par Notice that since $|\dot\gt|\geq |\dot\gt\bullet\gu|$, on account of \eqref{VarSPt} we clearly have $\F(\gt)\leq\Var_\SP(\gt)$, where the strict inequality holds in general, as $\dot\gt\bullet\gn\neq 0$ a.e. on $I_L$, when $\M$ has no ``flat" parts.
\br\label{RTheta} In Theorem~\ref{Tcomp}, we may and do assume that at each Jump point $s\in J_\Theta$, the Jump
$$ [\Theta]_s:=\Theta(s+)-\Theta(s-)$$
is bounded by $\pi$, i.e., $|[\Theta]_s|\leq\p$.
For this purpose, we consider the $\BV$ function $u=\E^{\IC \Theta}:I_L\to{\mathbb S}^1$ and build up an optimal lifting $\wid\Theta:I_L\to\gR$ of $u$ as in \cite{I}. Roughly speaking, we replace the Jump component $\Theta^J$ with a Jump function $\wid\Theta^J$ which has Jump set contained in $J_\Theta$ and such that for each $s\in J_\Theta$
$$ |[\wid\Theta^J]_s|\leq\p\,,\qquad [\wid\Theta^J]_s=[\Theta^J]_s+2k\pi\,,\quad k\in\gZ\,. $$
The optimal angle function is such that for a.e. $s\in I_L$ there exists $k\in\gZ$ such that $\wid\Theta(s)=\Theta(s)+2k\pi$, whence
$\cos\wid\Theta=\cos\Theta$ and $\sin\wid\Theta=\sin\Theta$ a.e. on $I_L$. This yields that formula \eqref{weakX} remains unchanged if we replace $\Theta$ with the optimal angle $\wid\Theta$. \er
%
%
%
%
\bt\label{Tangle} Under the hypotheses of Theorem~$\ref{Tcomp}$, and on account of Remark~$\ref{RTheta}$, we have
$$ |D\Theta| (I_L)=\F(\gt)\,.$$
More precisely, in the decomposition formula $|D\Theta| (I_L)=|D^a\Theta| (I_L)+|D^J\Theta|(I_L)+|D^C\Theta|(I_L)$ we have:
\beq\label{equalities} |D^a\Theta|(I_L)=\int_0^L|\dot\gt\bullet\gu|\,ds\,,\quad |D^C\Theta|(I_L)=|D^C\gt|(I_L)\,,\quad
|D^J\Theta|(I_L)=\sum_{s\in J_\gt}d_\SP(\gt(s+),\gt(s-))\,. \eeq
\et
\bpf
Let $\{P_h\}\sb\calP_\M(\gc)$ as in Theorem~\ref{Tcomp}, with transport vector fields $X_h$ and Darboux frames $(\gt_h,\gn_h,\gu_h)$, and let
$X$ be the $W^{1,1}$ transport vector field given by \eqref{weakX}.
%
%
\smallskip\par\noindent{\sc The a.c. components.} Recalling that $|D^a\Theta|(I_L)=\int_0^L|\dot\Theta|\,ds$, the first equality in \eqref{equalities} follows provided that
we show that for $\calL^1$-a.e. $s\in I_L$
\beq\label{eqAC}
 \dot\Theta(s)=\dot\gt(s)\bullet \gu(s)\,. \eeq
\par For this purpose, we first observe that from \eqref{weakX}, using that $X$ is a Sobolev function, and hence that it has a continuous representative, see eq. \eqref{DJX} below, it turns out that the Jump set of $\Theta$ agrees
with the Jump set of $\gt$ (and hence of $\gu$).
By the chain rule formula \eqref{chain} we infer that for a.e. $s$
$$ \dot X=\dot\Theta\,(-\sin\Theta\,\gt-\cos\Theta\,\gu)+\cos\Theta\,\dot\gt-\sin\Theta\,\dot\gu\,. $$
On the one hand, passing to the limit in the identities $\dot X_h\bullet\gt_h=0$ and $\dot X_h\bullet\gu_h=0$, by the a.e. convergences $X_h\to X$, $\gt_h\to\gt$, and $\gu_h\to\gu$, that hold true along subsequences, due to the $L^1$ convergences, we deduce that $\dot X\bullet\gt=0$ and $\dot X\bullet\gu=0$ a.e. on $I_L$. On the other hand, using that $|\gt_h|=1$, $|\gu_h|=1$, and $\gt_h\bullet\gu_h=0$, we also infer that $\dot\gt\bullet\gt=0$, $\dot\gu\bullet\gu=0$, and $\dot\gu\bullet\gt=-\dot\gt\bullet\gu$ a.e. on $I_L$. As in the proof of Proposition~\ref{PTheta}, by the above properties we obtain for a.e. $s$ the equations
$$ 0=\dot X\bullet \gt=-\sin\Theta\,(\dot\Theta-\dot\gt\bullet\gu)\,,\qquad 0=\dot X\bullet \gu=-\cos\Theta\,(\dot\Theta-\dot\gt\bullet\gu) $$
that clearly imply \eqref{eqAC}.
\smallskip\par\noindent{\sc The Cantor components.}
The second equality in \eqref{equalities} holds true if we show that
\beq\label{eqC} D^C\gt=\gu\,D^C\Theta\,. \eeq
\par To this aim, using again the chain rule formula \eqref{chain}, and since $X\in W^{1,1}$, we have
$$ 0=D^CX=-\sin\Theta\,\gt\,D^C\Theta-\cos\Theta\,\gu\,D^C\Theta+\cos\Theta\,D^C\gt-\sin\Theta\,D^C\gu$$
(where we choose good representatives of $\gt$, $\gu$, and $\Theta$) which is equivalent to the equation:
\beq\label{Cantor}
\cos\Theta\,(D^C\gt-\gu\,D^C\Theta)=\sin\Theta\,(D^C\gu+\gt\,D^C\Theta)\,. \eeq
Now, by taking the scalar products with $\gt$ and $\gu$ in equation \eqref{Cantor}, and observing that by \eqref{chain} we also have
$\gt\bullet D^C\gt=0$, $\gu\bullet D^C\gu=0$, and $\gt\bullet D^C\gu=-\gu\bullet D^C\gt$,
we obtain
$$ 0=\sin\Theta\,(-\gu\bullet D^C\gt +D^C\Theta)\,,\qquad \cos\Theta\,(\gu\bullet D^C\gt-D^C\Theta)=0 $$
which yields that $\gu\bullet D^C\gt=D^C\Theta$. But we have seen in Remark~\ref{RCantor} that $D^C\gt$ is tangential, namely, $D^C\gt=\gu\,(\gu\bullet D^C\gt)$.
Therefore, formula \eqref{eqC} is proved.
\smallskip\par\noindent{\sc The Jump components.} Recalling that $J_\gt=J_\gu=J_\Theta$, the third equality in \eqref{equalities} holds true if we show that for every $s\in J_\Theta$
\beq\label{eqJ}
|\Theta(s+)-\Theta(s-)|=d_\SP(\gt(s+),\gt(s-))\,. \eeq
%
%
\par Now, again by the chain rule formula \eqref{chain} we infer that
\beq\label{DJX} 0=D^J X=D^J(\cos\Theta\,\gt-\sin\Theta\,\gu)=\sum_{s\in J_\Theta}[\cos\Theta\,\gt-\sin\Theta\,\gu]_s\,\delta_s \eeq
where for each $s\in J_\Theta$
$$ [\cos\Theta\,\gt-\sin\Theta\,\gu]_s:= [\cos\Theta(s+)\,\gt(s+)-\sin\Theta(s+)\,\gu(s+)]-[ \cos\Theta(s-)\,\gt(s-)-\sin\Theta(s-)\,\gu(s-)]\,.  $$
For any fixed $s\in J_\Theta$, up to a rotation in the target space we may and do assume that $\gn(s)=(0,0,1)$, and hence we can write
$$ \gt(s\pm)=(\cos\a_\pm,\sin\a_\pm,0)\,,\qquad \gu(s\pm)=\gn(s)\tim \gt(s\pm)=(-\sin\a_\pm,\cos\a_\pm,0) $$
for some real numbers $\a_\pm$ satisfying $|\a_+-\a_-|\leq\p$.
Condition $[\cos\Theta\,\gt-\sin\Theta\,\gu]_s=0$ yields to the system
$$ \left\{\ba{l} \cos(\a_+-\Theta(s+))=\cos(\a_--\Theta(s-)) \\ \sin(\a_+-\Theta(s+))=\sin(\a_--\Theta(s-)) \ea \right. $$
which gives $\Theta(s+)-\Theta(s-)=\a_+ -\a_-$ mod $2\p$. By Remark~\ref{RTheta}, the optimal angle function satisfies $|\Theta(s+)-\Theta(s-)|\leq \p$. Since $|\a_+-\a_-|\leq\p$, we thus conclude that
$\Theta(s+)-\Theta(s-)=\a_+ -\a_-$. Therefore, equality \eqref{eqJ} follows by observing that $d_\SP(\gt(s+),\gt(s-))=|\a_+ -\a_-|$, as required.
\epf
\section{The high codimension case}\label{Sec:N}
In this section, we extend the previous results to the high codimension case of curves $\gc$ in $\M$, where $\M$ is a smooth (at least of class $C^3$), closed,
and compact immersed surface in $\RN$, with $N\geq 4$. We remark that $\M$ is not assumed to be oriented.
\par
We will only sketch the proofs: further details can be obtained by arguing in a way very similar to the codimension one case previously considered.
Moreover, when referring to analogous results from the previous sections, we shall tacitly assume that one has to replace $\SP$ and $\gR^3$ with $\SN$ and $\RN$, respectively, where $\SN$ is the unit hyper-sphere in $\RN$.
\adl\par\noindent{\large\sc Total curvature.} The Euclidean total curvature $\TC(\gc)$ of a curve $\gc$ in $\RN$ is defined as in the case $N=3$, and similar features hold.
Namely, if $\gc$ is smooth and regular, and $\gc:[0,L]\to\RN$
is its arc-length parameterization, one has $\TC(\gc)=\int_0^L|\gk|\,ds$,
where $\gk(s)$ is the curvature vector of $\gc$.
More generally, if $\gc$ has compact support and finite total curvature, then $\gc$ is rectifiable, and the tantrix $\gt=\dot\gc$ exists a.e., with
$\gt\in\SN$.
Moreover, the function $\gt:I_L\to\SN$ has bounded variation, and its essential variation in $\SN$ is equal to the total curvature of $\gc$, whereas formula \eqref{VarSPt} continues to hold for $\Var_{\SN}(\gt)$.
Furthermore, the total curvature of $\gc$ is equal to the limit of any sequence of polygonals $\{P_h\}$ in $\RN$ inscribed in $\gc$ and such that
$\mesh(P_h)\to 0$, see Remark~\ref{Rparall}.
\adl\par\noindent{\large\sc Total intrinsic curvature.}
If $\gc$ is a smooth and regular curve in $\M$, using that $\dot\gt\bullet\gt\equiv 0$, where $\bullet$ is the scalar product in $\RN$, the curvature vector $\gk(s):=\dot\gt(s)$ again decomposes as
\beq\label{curvN}
\gk(s)=\kk_g(s)\,\gu(s)+\kk_n(s)\,\gn(s)\,. \eeq
The unit conormal $\gu:[0,L]\to\SN$ is the unit vector orthogonal to $\gt$ and obtained by means of a positive rotation of $\gt$ on the tangent space $T_\gc\M$ along $\gc$, so that
$\gt\bullet\gu\equiv 0$ and the tangent space $T_{\gc(s)}\M$ is spanned by
$(\gt(s),\gu(s))$. Also, $\gn:[0,L]\to\SN$ is a smooth normal unit vector field (a section of the normal bundle).
\par The total intrinsic curvature $ \TC_{\M}(\gc)$ of a curve $\gc$ in $\M$ is defined as in the case $N=3$,
see Definition~\ref{Dcurv}.
For a polygonal $P$ in $\M$, we have $\TC_{\M}(P)=\gk_{\M}(P)$, and Proposition~\ref{Pappr} continues to hold.
\par Since inequality \eqref{curvnorm} is verified through the assumptions on $\M$, it turns out that a curve $\gc$ in $\M$ has finite total intrinsic curvature if and only if it has finite Euclidean total curvature, see \eqref{TCTCM}.
Whence, if $\TC_{\M}(\gc)<\ii$, then $\gc$ is rectifiable and one-sidedly smooth, see Remark~\ref{Rossm}, and the tantrix $\gt\in\BV(I_L,\SN)$.
\par As in the smooth case, we define the weak conormal $\gu$ in $\BV(I_L,\SN)$ by the unit vector orthogonal to $\gt$ and obtained by means of a positive rotation of $\gt$ on the tangent space $T_\gc\M$ along $\gc$.
\par Finally, formula $ D^C\gt=\gu\,(\gu\bullet D^C\gt)$ is obtained by arguing as in Remark~\ref{RCantor}, but this time observing that
$\gn_i\bullet D^C\gt=D^C(\gt\bullet\gn_i)=0$ for each $i=3,\ldots,N$, where $s\mapsto(\gn_3,\ldots,\gn_{N})(s)$ is a Lipschitz-continuous orthonormal frame that spans the normal space to $\M$ along $\gc$.
\adl\par\noindent{\large\sc Weak parallel transport.}
Proposition~\ref{PTheta} clearly extends to smooth and regular curves in $\M\sb\RN$.
Also, on account of \eqref{Theta} and \eqref{curvN}, by decomposing the derivative of the unit conormal
$$ \dot\gu=(\dot\gu\bullet\gt)\,\gt+\dot\gu^\perp$$
into the tangential and normal component to $\M$, and recalling that $\dot\gu\bullet\gt=-\gt\bullet\gu=-\dot\Theta$,
the parallel transport of (piecewise) smooth curves this time satisfies
\beq\label{dotXbis} \dot X= \cos\Theta\,\kk_n\,\gn-\sin\Theta\,\dot\gu^\perp\,, \eeq
where $\dot\gu^\perp=\dot\gu$ when $\gc$ is a geodesic arc.
\par Moreover, a compactness property as in Theorem~\ref{Tcomp} holds true: the limit function
$X\in W^{1,1}(I_L,\RN)$ satisfies \eqref{weakX}
for $\calL^1$-a.e. $s\in I_L$, where $\gt=\dot \gc$ is the unit tangent vector and
the conormal $\gu$ agrees with the weak-$^*$ $\BV$-limit of the sequence $\{\gu_h\}$ of the conormals to a subsequence of $\{P_h\}$.
\par In fact, compactness in $W^{1,1}$ is based on the validity of the estimate \eqref{traspest},
where the real constant $C_\M$ only depends on $\M$. Now, using this time the formula \eqref{dotXbis} for the derivative of $X=X_h$, where $\dot \gu_h^\perp=\dot\gu_h$, as $P_h$ is a polygonal in $\M$, the inequality
\eqref{traspest} is checked by arguing as in the proof of Theorem~\ref{Tcomp}, but this time observing that:
\ben \item the normal curvatures $\kk_n$ of the geodesics in $\M$, and their derivatives w.r.t. the arc-length parameter, are equibounded by a constant only depending on $\M$;
\item if $\gu$ is the unit conormal of a geodesic arc in $\M$, and $\dot\gu$ is its derivative w.r.t. the arc-length parameter, both $|\gu|$ and $|\dot\gu|$ are equibounded by a constant only depending on $\M$. \een
\br\label{Rest} The above properties follow from the smoothness and compactness of the surface $\M$ in $\RN$, and they will be discussed in
Sec.~\ref{Sec:Riem}, see also Example~\ref{Edotu}. \er
{\large\sc The angle function.} We now see that Theorem~\ref{Tangle} continues to hold.
For this purpose, by the structure \eqref{weakX} of the $W^{1,1}$ transport $X$ we again infer that $J_\Theta=J_\gt=J_\gu$.
Therefore, the equalities \eqref{equalities} hold true if we check the validity of the three formulas \eqref{eqAC}, \eqref{eqC}, and \eqref{eqJ}.
\par The equality \eqref{eqAC} involving the a.c. components is readily proved by means of the same argument.
\par As to the Cantor components, using that $\gt\bullet D^C\gt=0$, $\gu\bullet D^C\gu=0$, and
$\gt\bullet D^C\gu=-\gu\bullet D^C\gt$, we similarly obtain that $\gu\bullet D^C\gt=D^C\Theta$, whence the equality \eqref{eqC} follows since we have already checked
the tangential property $D^C\gt=\gu\,(\gu\bullet D^C\gt)$.
\par As to the Jump components, for any $s\in J_\Theta$, up to a rotation we may and do assume that the tangent space $T_{\gc(s)}\M$ is spanned by the first two vectors of the canonical basis in $\RN$. Therefore, we can write
$$ \gt(s\pm)=(\cos\a_\pm,\sin\a_\pm,0_{\gR^{N-2}})\,,\qquad \gu(s\pm)=(\sin\a_\pm,-\cos\a_\pm,0_{\gR^{N-2}}) $$
for some real numbers $\a_\pm$ satisfying $|\a_+-\a_-|\leq\p$, whence $d_\SN(\gt(s+),\gt(s-))=|\a_+-\a_-|\leq\pi$.
Condition $[\cos\Theta\,\gt-\sin\Theta\,\gu]_s=0$ implies again that $\Theta(s+)-\Theta(s-)=\a_- -\a_+$ mod $2\p$, whereas Remark~\ref{RTheta} on the optimal angle function $\Theta$ continues to hold, whence equation \eqref{eqJ} is satisfied and the proof is complete.
%
%
\section{Gauss-Bonnet theorem and representation formula}\label{Sec:GB}
In this section, we discuss the validity of Gauss-Bonnet formula in the setting of domains in $\M$ bounded by simple and closed curves with finite total curvature, Theorem~\ref{TGB}.
As a consequence, we shall obtain an explicit representation formula for the total intrinsic curvature of curves in immersed surfaces, Theorem~\ref{Tequality}.
\adl\par\noindent{\large\sc The Gauss-Bonnet theorem.} We have:
\bt\label{TGB} Let $\M$ be a smooth, closed, compact, and immersed surface in $\RN$, where $N\geq 3$.
Let $\gc:[0,L]\to\M$ be a simple and closed rectifiable curve with finite total curvature, $\TC_\M(\gc)<\ii$.
Let $k(s)\,ds:=D\Theta[0,s)$, where $\Theta$ is the left-continuous representative of the optimal angle function of the parallel transport along $\gc$, see Theorems~$\ref{Tcomp}$ and $\ref{Tangle}$, so that
$$ \int_0^L k(s)\,ds=\Theta(L)-\Theta(0)\,. $$
Let $U$ be the open set in $\M$ enclosed by the oriented curve $\gc$. Moreover, assume that $U$ is simply connected, and that for a.e. $s\in I_L$ the tangent vector $\gt(s)$ is positively oriented w.r.t. the natural orientation on the boundary of $U$ at $\gc(s)$.
Finally, let $\gK$ denote the Gauss curvature of $\M$, and $\a$ the oriented angle from $\gt(L-)$ to $\gt(0+)$ at the junction point $\gc(0)=\gc(L)$.
Then we have:
$$ \int_U\gK\,dA=2\p-\int_0^L k(s)\,ds-\a \,.$$
\et
\par Notice that if $\gc$ is smooth, by Proposition~\ref{PTheta} we know that $D\Theta=\dot\Theta\,\calL^1$, with $\dot\Theta(s)=\kk_g(s)$ for each $s$,
so that we recover the classical formula, as $\int_0^L k(s)\,ds=\int_\gc\,\kk_g(s)\,ds$, see \eqref{GB}. In a similar way one may proceed in the case of piecewise smooth curves, this
time obtaining an extra term given by the sum of the oriented turning angles at the corner points of $\gc$, in correspondence to the Jump points of the angle function $\Theta$ in $I_L$, plus a possible extra term at the junction point $\gc(0)=\gc(L)$.
Therefore, our Theorem~\ref{TGB} extends the classical Gauss-Bonnet theorem to the wider class of curves with finite total curvature.
\par
If $\TC_\M(\gc)=\ii$, in fact, we expect that there is no way to find a finite measure that contains the information
(given by the derivative $D\Theta$ of the angle function of the parallel transport along the curve) on the ``signed geodesic curvature" of the curve $\gc$.
\par Finally, a more general result could be obtained if $U$ fails to be simply-connected, assuming $\M$ oriented. This time, the term $2\pi\cdot\chi(U)$ appears,
$\chi(U)$ being the Euler-Poincar\'e characteristic of $U$.
\adl\par\noindent\bpff{\sc of Theorem~\ref{TGB}:}
Let $\{P_h\}\sb\calP_\M(\gc)$ as in Theorem~\ref{Tcomp}, with transport vector fields $X_h:[0,L]\to\M$ given by \eqref{Xh}.
Let $U_h$ be the open set in $\M$ enclosed by the oriented closed polygonal $P_h$, and $\gi_h(x)$ the index of $P_h$ at the point $x\in\M$.
By uniform convergence, for $h$ sufficiently large we can choose a simply-connected and open set $U_h$ in $\M$ such that the index $\gi_h$ is equal to zero outside $U_h$.
By applying the classical Gauss-Bonnet theorem, and recalling that by our assumptions $P_h(0)=P_h(L)=\gc(0)=\gc(L)$, it is readily checked that the equality
$$ \int_{U_h} \gi_h\,\gK\,dA=2\p-\int_0^L k_h(s)\,ds-\a_h $$
holds true, where $k_h(s):=D\Theta_h[0,s)$, so that
$\int_0^L k_h(s)\,ds=\Theta_h(L)-\Theta_h(0)$, and $\a_h$ is the oriented angle from $\gt_h(L)$ to $\gt_h(0)$ at the junction point $P_h(0)=P_h(L)$.
By the weak-$^\ast$ convergence of $D\Theta_h$ to $D\Theta$, we infer that
$\int_0^L k_h(s)\,ds\to \int_0^L k(s)\,ds$ as $h\to\ii$. On the other hand, by the uniform convergence of $P_h$ to $\gc$ we obtain that $\int_{U_h}\gi_h\,\gK\,dA\to \int_{U}\gK\,dA$. Finally, since $\gc$ is one-sidedly smooth, we also infer that $\a_h\to\a$, as required.
\epff
%
{\large\sc The representation formula.}
In general, by the sequential lower-semicontinuity of the total variation w.r.t. the weak-$^\ast$ convergence, in Theorem~\ref{Tcomp}
(that holds true for curves contained in surfaces $\M$ of $\RN$) we only have
$$
|D\Theta| (I_L)\leq \lim_{h\to\ii}|D\Theta_h| (I_L)=\lim_{h\to\ii}\gk_\M(P_h)=\TC_\M(\gc) $$
where the last equality follows from Proposition~\ref{Pappr}.
\par
As a consequence, by Theorem~\ref{Tangle} we obtain the inequality
\beq\label{TMgeqF} \TC_\M(\gc)\geq\F(\gt) \eeq
where $\F(\gt)$ is the energy functional given by \eqref{Ft},
and we expect that equality holds in \eqref{TMgeqF} in full generality.
%
\par In fact, for piecewise smooth and regular curves $\gc$ in $\M$, one has:
$$ \F(\gt)=\int_0^L|\kk_g(s)|\,ds+\sum_i|\a_i| $$
so that it suffices to apply Theorem~\ref{Tdens} and \eqref{TCMsmooth}.
\br\label{Rmon} We now readily check that {\em equality holds in \eqref{TMgeqF} for convex or concave curves} with finite total intrinsic curvature, i.e., for simple and closed curves $\gc$ such that the
right-hand (or left-end) side region with boundary the trace of $\gc$ is a geodesically-convex subset of $\M$.
For non-closed curves, this means that all the length minimizing arcs connecting two points of the curve lie on the same side w.r.t. the tantrix of the curve.
\par In this case, in fact,
for any polygonal $P_h$ in $\M$ inscribed in $\gc$, the angle $\Theta_h$ of the parallel transport along $P_h$
is a monotone function. Therefore, for each $(a,b)\sb I_L$ we have $|D\Theta_h|(a,b)=|\Theta_h(b-)-\Theta_h(a+)|$.
The a.e. convergence of $\Theta_h$ to $\Theta$, that holds true for a subsequence, yields that the angle $\Theta$ is a monotone function, too,
whence $|D\Theta|(a,b)=|\Theta(b-)-\Theta(a+)|$. As a consequence, we obtain the strict convergence $|D\Theta_h|(I)\to|D\Theta|(I)$, which implies
the equality sign in \eqref{TMgeqF}, on account of Theorem~\ref{Tangle}. \er
\par By exploiting (in Proposition~\ref{PGB}) the generalized Gauss-Bonnet theorem~\ref{TGB}, we are able to prove that equality holds in \eqref{TMgeqF},
even in the non trivial case of surfaces $\M$ with positive Gauss curvature.
%
\bt\label{Tequality} Let $\M$ be a smooth (at least of class $C^3$), closed, and compact (not necessarily oriented) immersed surface in $\RN$.
Then, for every rectifiable curve $\gc$ in $\M$ with finite total curvature, $\TC_\M(\gc)<\ii$, we have
$$ \TC_\M(\gc)=\F(\gt) $$
where $\F(\gt)$ is given by \eqref{Ft} and $\gt=\dot\gc$ is the tantrix of the curve.
\et
\par We first observe that Theorem~\ref{Tequality} holds true as a consequence of the following proposition, that will be proved in the second part of this section.
\bp\label{PGB} Let $\gc:[0,L]\to\M$ be a rectifiable curve with finite total curvature (parameterized by arc-length), and let $\Theta$ denote the left-continuous representative of the optimal angle of the parallel transport $X$ along $\gc$, with initial condition $X(0)=\gt(0)$.
Let
$\{P_h\}\sb\calP_\M(\gc)$ with modulus $\m_\gc(P_h)\to 0$.
Assume that $P_h$ is generated by the consecutive vertexes $\gc(s_i)$, where $0=s_0<s_1<\cdots < s_n=L$ (with $\{s_i\}$ and $n$ depending on $h$), and that every $s_i$ is not a Jump point of the angle function $\Theta$.
Also, let $\Theta_h$ denote the angle of the parallel transport $X_h$ along $P_h$, with initial condition $X_h(0)=\gt(0)$.
Then, for $h$ sufficiently large there exists a piecewise constant function $\wid\Theta_h:I_L\to\gR$ such that:
\bi
\item[$(a)$] for each $i=1,\ldots,n$, there exists a parameter $\wid s_i\in[s_{i-1},s_i[$ such that $\wid\Theta_h(s)=t_i\,\Theta(\wid s_i+)+(1-t_i)\Theta(\wid s_i-)$ for any $s\in]s_{i-1},s_i[$, where $t_i\in[0,1]$\,;
\item[$(b)$] $\Var(\Theta_h)\leq\Var(\wid\Theta_h)+\e_h$, where $\e_h\to 0^+$ as $h\to\ii$.
\ei \ep
%
\bpff{\sc of Theorem~\ref{Tequality}:} We first notice that the assumption on the continuity of the angle function $\Theta$ at the points $s_i$
is assumed without loss of generality, as the Jump set $J_\Theta$ is at most countable.
\par
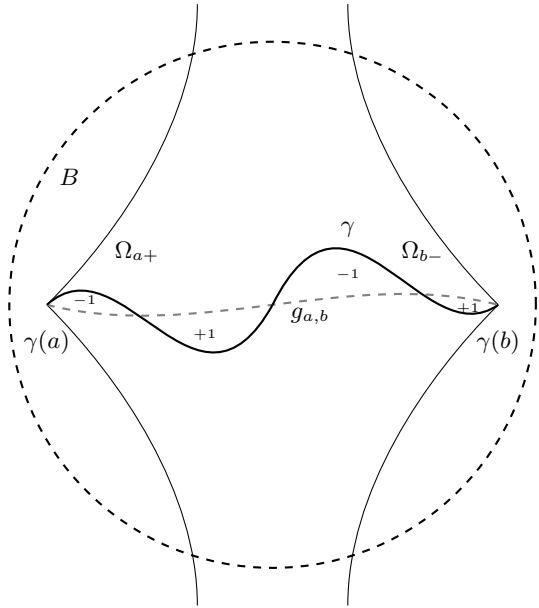
\begin{figure}[ht]
	\centering
		\begin{tikzpicture}

\draw (0,0) .. controls (1,1) and (2,2.5) .. (2,4);
\draw (0,0) .. controls (1,-1) and (2,-2.5) .. (2,-4);
\draw[thick,dashed] (3,0) circle (3.5);
\node at (0.3,1.7) {\small $B$};
\node at (0,-0.5) {\small $\gamma(a)$};
\node at (6,-0.5) {\small $\gamma(b)$};
\node at (1.2,0.7) {\small $\Omega_{a+}$};
\node at (5,0.7) {\small $\Omega_{b-}$};
\draw (6,0) .. controls (5,1) and (4,2.5) .. (4,4);
\node at (4,1) {\small $\gamma$};
\draw (6,0) .. controls (5,-1) and (4,-2.5) .. (4,-4);
\draw[thick,gray,dashed] (0,0) .. controls (2,-0.5) and (4,0.5) .. (6,0);
\node at (3.5,-0.15) {\small $g_{a,b}$};
\draw[thick] (0,0) .. controls (1,0.9) and (2,-1.8) .. (3,0);
\draw[thick] (6,0) .. controls (5,-0.7) and (4,2) .. (3,0);
\node at (0.5,0.05) {\tiny $-1$};
\node at (5.6,-0.05) {\tiny $+1$};
\node at (2.1,-0.4) {\tiny $+1$};
\node at (4,0.4) {\tiny $-1$};
\end{tikzpicture}
	\caption{The simply-connected closed set $\Omega=B\cap\Omega_{a+}\cap\Omega_{b-}$ of Lemma~\ref{LGB}. The arc $\g$ is drawn with a continuous line, and the geodesic arc connecting $\g(a)$ and $\g(b)$
with a dashed line.}
	\label{GBpic}
\end{figure}

 Property~$(a)$ in Proposition~\ref{PGB} implies that the modified angle $\wid\Theta_h$ is a competitor to the computation of the essential variation of $\Theta$, compare \cite[Sec.~3.2]{AFP}, whence
$\Var(\wid\Theta_h)\leq\Var(\Theta)$. By property~$(b)$ in Proposition~\ref{PGB}, we deduce that $\limsup_h\Var(\Theta_h)\leq\Var(\Theta)$. The weak convergence of $\Theta_h$ to $\Theta$, see Theorem~\ref{Tcomp}, yields that $\Var(\Theta)\leq\liminf_h\Var(\Theta_h)$, whence we obtain the strict convergence $\Var(\Theta_h)\to\Var(\Theta)$.
Since $\Var(\Theta_h)=\gk_\M(P_h)$, whereas by Proposition~\ref{Pappr} we know that $\gk_\M(P_h)\to\TC_\M(\gc)$, and by Theorem~\ref{Tangle} that $\Var(\Theta)=|D\Theta|(I_L)=\F(\gt)$, we conclude that $\TC_\M(\gc)=\F(\gt)$, as required. \epff
{\large\sc A localization lemma.} Proposition~\ref{PGB} will be proved by exploiting Theorem~\ref{TGB}, see formulas \eqref{GBlocal} and \eqref{GBglobal}.
For this purpose, we shall make use of the following result, which is illustrated in Figure~$\ref{GBpic}$.
\bl\label{LGB} Given any one-sidedly smooth curve $\gamma:[0,L]\to\M$, parameterized in arc length, there is $\varepsilon_0>0$ such that for any $[a,b]\subset[0,L]$
satisfying $b-a<\varepsilon_0$ we can find a simply-connected closed set $\Omega\subset\M$ for which $\gamma([a,b])\subset\Omega$ and
$\gamma(a),\gamma(b)\in\pa\Omega$, in such a way that the minimal geodesic arcs connecting any couple of points in the curve $\gamma([a,b])$ are contained in $\Omega$. In particular, the geodesic arc connecting $\gamma(a)$ and $\gamma(b)$ divides $\Omega$ in two connected components. \el
\bpff{\sc of Lemma~\ref{LGB}:} Let us fix $s\in[0,L]$. Let $\varepsilon_1(s)$ be half of the injectivity radius of $\M$ at $\gamma(s)$, and
let $\varepsilon_1 :=\ \inf_{s\in[0,L]} \varepsilon_1(s)$,
so that by compactness of the curve $\g$ and smoothness of $\M$, which implies a uniform bound on the sectional curvature of $\M$, we get 
$$\varepsilon_1\ =\ \min_{s\in[0,L]} \varepsilon_1(s)\,,\qquad \e_1>0\,.$$

From now on we will consider only points $a,b\in[0,L]$ at most $\varepsilon_1$ apart, so that the geodesic $g_{a,b}$ from $\gamma(a)$ to $\gamma(b)$ can be uniquely defined as the shortest path connecting $\gamma(a)$ and $\gamma(b)$.

We first extend $\gamma$ to a one-sidedly smooth curve defined in a neighborhood of $[0,L]$.
Now, if $a\in[0,L)$, the right geodesic tangent at $\gamma(a)$, i.e., the geodesic $g_{a+}$ starting in $\gamma(a)$ with tangent vector $\gt(a+)$, is well-defined. Moreover, by Remark~\ref{Rossm}, where we fix e.g. $\delta=\p/4$, and by the smoothness and compactness on $\M$, it turns out that for any $s\in[0,L]$, there is $\e_2(s)\in(0,\varepsilon_1]$ such that if $a,b\in[s,s+\e_2(s)]$, with $a<b$, then the angle in $\gamma(a)$ between $g_{a+}$ and $g_{a,b}$ is less than $\delta/2$. Let
$$\varepsilon_2 \ :=\ \inf_{s\in[0,L]} \varepsilon_2(s)\,,$$
so that $\e_2$ is a positive minimum, $\e_2>0$, by continuity of the function $\e_2(s)$ in $[0,L]$.

By the previous construction, if $0\leq a<c<b\leq L$ are chosen so that $b-a\leq\varepsilon_2$, then the angle between $g_{a,b}$ and $g_{a,c}$ in $\gamma(a)$ is smaller than $\delta$. As a consequence, the curve $\gamma([a,b])$ is contained in the geodesic sector $\Omega_{a+}$ bounded by the geodesics from $\gamma(a)$ with starting direction tilted by $\pm\delta$ from the one of $g_{a,b}$.

With the same reasoning applied to $b\in(0,L]$ and to the left geodesic $g_{b-}$, we can find a positive number $\varepsilon_3\in(0,\varepsilon_1]$ such that if $0\leq a<c<b\leq L$ satisfy $b-a\leq\varepsilon_3$, then the angle between $g_{a,b}$ and $g_{c,b}$ in $\gamma(b)$ is smaller than $\delta$. Hence the curve $\gamma([a,b])$ is contained in the geodesic sector $\Omega_{b-}$ bounded by the geodesics from $\gamma(b)$ with starting direction tilted by $\pm\delta$ from the one of $g_{a,b}$.

Let then $\varepsilon_0:=\min\{\varepsilon_2,\varepsilon_3\}$, and let $\Omega:=\Omega_{a+}\cap\Omega_{b-}\cap B$, where $B$ is the intersection of the geodesic balls of radii $\varepsilon_0$ centered in $\gamma(a)$ and $\gamma(b)$, see Figure~\ref{GBpic}.
We thus conclude that if $a,b\in[0,L]$ are such that $0<b-a<\varepsilon_0$, then $\gamma([a,b])\subset\Omega$, the closed set $\Omega$ is simply-connected, and $\gamma(a),\gamma(b)\in \pa\Omega$.
Moreover, the minimal geodesic arcs connecting any couple of points in the curve $\gamma([a,b])$ are contained in $\Omega$.
 Finally, the arc $g_{a,b}$ divides $\Omega$ in two connected components, as required.
\epff
{\large\sc The role of Gauss-Bonnet theorem.} In order to make the proof of Proposition~\ref{PGB} more clear, we first recall how the equality $\TC(\gc)=\Var_{\SN}(\gt)$ is checked for curves $\gc$ in $\RN$ with finite total curvature, and then
deal with the case $N=2$, where we apply a ``planar" version of the Gauss-Bonnet theorem~\ref{TGB}.
\br\label{Rparall} Let $P_h$ be an inscribed polygonal to the curve $\gc:[0,L]\to\RN$ (parameterized by arc-length) and generated by the consecutive vertexes $\gc(s_i)$, where $0=s_0<s_1<\cdots < s_n=L$, and let $\gv_i$ be the oriented segment of $P_h$ from $\gc(s_{i-1})$ to $\gc(s_i)$. If $\gt_h$ is the tantrix of $P_h$ in $\SN$, the value of $\gt_h$ in $\gv_i$ is an average of the values of the restriction of the tantrix $\gt$ of $\gc$ to $(s_{i-1}, s_i)$, when completed to a continuous curve in $\SN$ by connecting with geodesic arcs the points $\gt(s-)$ and $\gt(s+)$ for each $s\in J_\gt\cap(s_{i-1},s_i)$, compare \cite{AR}.
This property implies that $\Var_\SN(\gt_h)\leq\Var_\SN(\gt)$. If $\{P_h\}$ is an inscribed sequence satisfying $\mesh(P_h)\to 0$, the weak $\BV$ convergence of $\gt_h$ to $\gt$ implies the lower semicontinuity inequality $\Var_\SN(\gt)\leq\liminf_h\Var_\SN(\gt_h)$, yielding the strict convergence $\Var_\SN(\gt_h)\to\Var_\SN(\gt)$. Using that $\Var_\SN(\gt_h)\to\TC(\gc)$, one concludes that $\TC(\gc)=\Var_{\SN}(\gt)$. \er
%
\par When $\gc$ is a {\em planar curve}, i.e., when $N=2$, the value of $\gt_h\in\Su$ on the segment $\gv_i$
is equal to one of the values of the ``completion" in $\Su$ of the restriction of the tantrix $\gt$ to the interval $]s_{i-1}, s_i[$.
\par We now see that this property can be rewritten in terms of angle functions, and hence of the ``planar" version of the Gauss-Bonnet
theorem~\ref{TGB}, where of course $\gK\equiv 0$.
This is the starting point to treat the case of curves on surfaces. In the proof of Proposition~\ref{PGB}, moreover, we have to consider the angle of the
parallel transport, and to deal with the extra term given by the integral of the Gauss curvature.
\par We thus denote by $\om(s)$ the oriented angle from $\gt(s)$ to the fixed direction $\gt(0)$, where we choose $\gt$ equal to the left-continuous representative of the $\BV$-function $\dot\gc$. We assume moreover that $P_h:[0,L]\to\gR^2$ is parameterized with constant velocity on each interval $]s_{i-1},s_i[$, in such a way that
$P_h(s_{i})=\gc(s_{i})$ for each $i$, and that
every $s_i$ is not a Jump point of $\gt$.
\par
If $\om_h(s)$ is the oriented angle from $\gt_h(s)$ to $\gt(0)$, then $\om_h(s)$ is constant on each interval $]s_{i-1},s_i[$.
In order to show that $\Var(\om_h)\to\Var(\om)$, by \cite[Lemma~1]{Re} we may and do assume that $\gc$ is a simple arc.
Also, by Lemma~\ref{LGB} we can reduce to the following situation, for $h$ large enough.
\par Denote by $\angle{\gt(s)\gv_i}$ the oriented angle from $\gt(s)$ to $\gv_i$, where $s\in[s_{i-1},s_i[$, and $\gv_i$ is the oriented segment of $P_h$ from $\gc(s_{i-1})$ to $\gc(s_i)$.
For $i=1,\ldots,n$, letting $\a_i:=\angle{\gt(s_{i-1})\gv_i}$, if $\a_i\neq 0$, we choose the {\em first} parameter
$\ol s_i$ in the interval $]s_{i-1},s_i]$ such that $\gc(\ol s_i)\in \gv_i$.
Then, by Lemma~\ref{LGB}, {\em the angle $\ol\be_i:=\angle{\gt(\ol s_i)\gv_i}$ cannot have the same sign as $\a_i$}, i.e., $\a_i\cdot\ol\be_i\leq 0$.
Moreover, denoting by $\g_i$ the oriented closed curve given by the join of the arc $\gc_i:=\gc_{\vert[s_{i-1},\ol s_i]}$ plus the segment of $P_h$ from $\gc(\ol s_i)$
to $\gc(s_{i-1})$, {\em the index of $\g_i$ on the open set $U_i$ enclosed by $\g_i$ is equal to the sign of $\a_i$}, see Figure~$\ref{GBpic}$.
We thus have
$$ \om(\ol s_{i})-\om(s_{i-1})=\a_i-\ol\be_i\,,\qquad \a_i\neq 0\,,\quad\a_i\cdot\ol\be_i\leq 0\,.$$
\par Letting now $f_i(s):=\om(s)-\om(s_{i-1})$, we get $f_i(s_{i-1})<\a_i$ and $f_i(\ol s_{i})\geq\a_i$, when $\a_i>0$ and $\ol\be_i\leq0$, whereas
$f_i(s_{i-1})>\a_i$ and $f_i(\ol s_{i})\leq\a_i$, when $\a_i<0$ and $\ol\be_i\geq0$.
Therefore, using that $\om$ is a function with bounded variation,
we find $\wid s_i\in]s_{i-1},\ol s_i[$ such that either
$\a_i=t_i\,f_i(\wid s_i+)+(1-t_i)f_i(\wid s_i-)$ for some $t_i\in[0,1]$, if $\wid s_i$ is a Jump point of $f_i$, or $\a_i=f_i(\wid s_i)$, otherwise. When $\a_i=0$, we clearly have $\a_i=f_i(0)$.
\par Recall that $\om(s_0)=0$ and $\a_i:=\angle{\gt(s_{i-1})\gv_i}$. Setting $\be_i:=\angle{\gt(s_i)\gv_i}$, by the previous discussion based on Lemma~\ref{LGB}, we also get:
$$ \om(s_j)-\om(s_{j-1})=\a_j-\be_j\qquad\fa\,j=1,\ldots,n\,. $$
Moreover, for $j=1,\ldots,n-1$, the oriented turning angle of the polygonal $P_h$ at the corner point $\gc(s_j)$ is equal to $\a_{j+1}-\be_j$.
We thus have $\om_h(s)=\a_1$ if $s\in]s_0,s_1[$, whereas if $s\in]s_{i-1},s_i[$, and $i=2,\ldots,n$, then
$$\om_h(s)= \a_1+\sum_{j=1}^{i-1}(\a_{j+1}-\be_j)=\a_i+\sum_{j=1}^{i-1} (\a_j-\be_j)=\a_i+\sum_{j=1}^{i-1} (\om(s_j)-\om(s_{j-1}))=\a_i+\om(s_{i-1})\,. $$
We thus conclude that for each $i=1,\ldots,n$ there exists $\wid s_i\in[s_{i-1},s_i[$ and $t_i\in[0,1]$ such that
$$\om_h(s)=t_i\,\om(\wid s_i+)+(1-t_i)\,\om(\wid s_i-) \qquad \fa\,s\in]s_{i-1},s_i[\,. $$
\par The above property, that actually expresses the parallelism condition in term of angle functions, implies that $\om_h$ is a competitor to the computation of the essential variation of $\om$, whence $\Var(\om_h)\leq\Var(\om)$.
By the weak-$^*$ $\BV$ convergence of $\om_h$ to $\om$, which ensures that $\Var(\om)\leq\liminf_h\Var(\om_h)$, we obtain the strict convergence $\Var(\om_h)\to\Var(\om)$.
\adl\par\noindent\bpff{\sc of Proposition~\ref{PGB}:} By \cite[Lemma~1]{Re}, the curve $\gc$ being one-sidedly smooth, it consists of finitely many simple arcs. Therefore, we clearly may and do assume that $\gc$ is a simple arc.
\par We let $P_h:[0,L]\to\M$ be parameterized with constant velocity on each interval $]s_{i-1},s_i[$, in such a way that $P_h(s_{i})=\gc(s_{i})$ for each $i$.
Notice that by the uniform convergence of $P_h$ to $\gc$, for $h$
sufficiently large the subset of $\M$ enclosed by the curves $\gc$ and $P_h$ is a simply-connected domain $U_h$ of $\M$ with small surface area.
In particular, $U_h$ can be equipped with an orientation, that is inherited by the tangent space $T_{\gc(s)} \M$ along the curve.
If $\gv_0,\gv_1\in T_{\gc(s)}\M$ are non-trivial vectors, we shall thus denote by $\angle{\gv_0\gv_1}$ the oriented angle in $T_{\gc(s)}\M$ from $\gv_0$ to $\gv_1$, for any $s\in[0,L]$.
The rest of the proof is divided into three steps.
\smallskip\par\noindent{\sc Step~1:} We prove property~(a) in Proposition~\ref{PGB}.
\par Choose $h$ large enough so that $\m_\gc(P_h)\leq\e_0$, where the positive constant $\e_0>0$ is given by Lemma~\ref{LGB} in correspondence to the curve $\gc$.
We are now in a situation similar to the one described in the planar case.
\par For $i=1,\ldots,n$, letting $\a_i:=\angle{\gt(s_{i-1})\gt_h(s_{i-1}+)}$, if $\a_i\neq 0$, we choose the {\em first} parameter
$\ol s_i$ in the interval $]s_{i-1},s_i]$ such that $\gc(\ol s_i)=P_h(\wih s_i)$ for some $\wih s_i\in ]s_{i-1},s_i]$, and let
$\ol\be_i:=\angle{\gt(\ol s_i)\gt_h(\wih s_i-)}$.
\par By Lemma~\ref{LGB}, the angle $\ol\be_i$ cannot have the same sign as $\a_i$, i.e., $\a_i\cdot\ol\be_i\leq 0$, see Figure~\ref{GBpic}.
Also, denoting by $\wid\g_i$ the oriented closed curve given by the join of the arc $\gc_i:=\gc_{\vert[s_{i-1},\ol s_i]}$ plus the geodesic arc of $P_h$ reversely oriented from $\gc(\ol s_i)$
to $\gc(s_{i-1})$, the index of $\wid\g_i$ on the open set $\wid U_i$ enclosed by $\wid\g_i$ is equal to $\pm 1$, in concordance with the sign of the initial angle $\a_i$,
see Figure~$\ref{GBpic}$.
\par Therefore, the Gauss-Bonnet theorem~\ref{TGB} yields:
\beq\label{GBlocal} \left\{\ba{ll} \ds\Theta(\ol s_{i})-\Theta(s_{i-1})=\a_i-\ol\be_i-\int_{\wid U_i}\gK\,dA &\quad \If\quad \a_i>0 \\
\ds\Theta(\ol s_{i})-\Theta(s_{i-1})=\a_i-\ol\be_i+\int_{\wid U_i}\gK\,dA &\quad \If\quad \a_i<0 \ea\right. \eeq
where, we recall, $\a_i\cdot\ol\be_i\leq 0$.
\par We first consider the easier case when $\gK\leq 0$.
Letting $\ol\a_i:=\a_i-\int_{\wid U_i}\gK\,dA$, if $\a_i>0$, and
$\ol\a_i:=\a_i+\int_{\wid U_i}\gK\,dA$, if $\a_i<0$, in both cases the sign of $\ol\a_i$ is concordant with the sign of $\a_i$, and definitely:
$$ \Theta(\ol s_{i})-\Theta(s_{i-1})=\ol\a_i-\ol\be_i\,,\qquad \ol\a_i\neq 0\,,\quad \ol\a_i\cdot\ol\be_i\leq 0 \,. $$
Denoting $f_i(s):=\Theta(s)-\Theta(s_{i-1})$, we get $f_i(s_{i-1})<\ol\a_i$ and $f_i(\ol s_{i})\geq\ol\a_i$, when $\ol\a_i>0$ and $\ol\be_i\leq0$, whereas
$f_i(s_{i-1})>\ol\a_i$ and $f_i(\ol s_{i})\leq\ol\a_i$, when $\ol\a_i<0$ and $\ol\be_i\geq0$.
Therefore, recalling that the angle function $\Theta$ has bounded variation,
and setting $\Theta_{h,i}:=\ol\a_i+\Theta(s_{i-1})$, in both cases we find $\wid s_i\in]s_{i-1},\ol s_i[$ such that
$$\Theta_{h,i}=t_i\,\Theta(\wid s_i+)+(1-t_i)\,\Theta(\wid s_i-) $$
for some $t_i\in[0,1]$, if $\wid s_i$ is a Jump point of $\Theta$, or $\Theta_{h,i}=\Theta(\wid s_i)$, otherwise.
When $\a_i=0$, we clearly have $\a_i=f_i(0)$, and we obviously choose $\Theta_{h,i}:=\Theta(s_{i-1})$.
\par In order to treat the general case, where the Gauss curvature $\gK$ may possibly take positive values, in Step~2 we shall prove the following:
\smallskip\par\noindent{\bf Claim.} {\em For each $i=1,\ldots,n$, we can find a coefficient $\lambda_i\in[-1,1]$ such that with
$\ol\a_i:=\a_i+\lambda_i\int_{\wid U_i}\gK\,dA$ and $\Theta_{h,i}:=\ol\a_i+\Theta(s_{i-1})$, we have
$$\Theta_{h,i}=t_i\,\Theta(\wid s_i+)+(1-t_i)\,\Theta(\wid s_i-)$$
for some $\wid s_i\in[s_{i-1}, s_i[$ and $t_i\in[0,1]$.}
\smallskip\par
Setting in fact
\beq\label{widT} \wid\Theta_h(s):=\Theta_{h,i}\qquad\If\quad s\in]s_{i-1},s_i[\,,\qquad\fa\,i=1,\ldots,n \eeq
property~(a) in Proposition~\ref{PGB} holds true.
\smallskip\par\noindent{\sc Step~2:} We prove the Claim, by generalizing the previous argument. Denote for simplicity
$$\Delta\Theta_i:= \Theta(\ol s_{i})-\Theta(s_{i-1})\,,\qquad K_i:=\int_{\wid U_i}\gK\,dA\,.  $$
If $K_i\leq0$, we argue exactly as in Step~1, so that we now assume $K_i>0$.
\par We first consider the case $\a_i>0$ and $\ol\be_i\leq 0$, so that the first equation in \eqref{GBlocal} becomes
\beq\label{GBlocalbis} \Delta\Theta_i=\a_i-K_i-\ol\be_i \eeq
and we can write $\a_i=\lambda\,K_i$ for some $\lambda>0$.
We now distinguish among the possible values of the term $\Delta\Theta_i$.
\ben\item If $\Delta\Theta_i=0$, then by \eqref{GBlocalbis} we get $\lambda\in]0,1]$ and $\ol\be_i=(\lambda-1)\,K_i$.
Letting $\ol\a_i:=\a_i-\lambda\,K_i$, we clearly have $\ol\a_i=0=f_i(0)$, where $f_i(s)$ is defined as in Step~1.
\item If $\Delta\Theta_i>0$, then $\ol\be_i=-\m\,K_i$ for some $\m\geq 0$, so that \eqref{GBlocalbis} becomes $\Delta\Theta_i=(\lambda+\m-1)\,K_i$, whence $\lambda+\m>1$.
If $\lambda\geq 1$, letting $\ol\a_i:=\a_i-K_i$, we have
$$f_i(0)\leq\ol\a_i\,,\qquad f_i(s_{i-1})=\Delta\Theta_i=\ol\a_i+\m\,K_i\geq\ol\a_i \,. $$
If $\lambda\in]0,1]$, instead, letting $\ol\a_i:=\a_i-\lambda\,K_i$ we again have $\ol\a_i=0=f_i(0)$.
\item If $\Delta\Theta_i<0$, by \eqref{GBlocalbis} we have $\a_i-\ol\be_i<K_i$, hence $\lambda\in[0,1[$, so that we again let $\ol\a_i:=\a_i-\lambda\,K_i=f_i(0)$. \een
\par We now deal with the case $\a_i<0$ and $\ol\be_i\geq 0$, so that the second equation in \eqref{GBlocal} becomes
\beq\label{GBlocalter} \Delta\Theta_i=\a_i+K_i-\ol\be_i \eeq
and hence this time $\a_i=-\lambda\,K_i$ for some $\lambda>0$.
\ben\item If $\Delta\Theta_i=0$, then by \eqref{GBlocalter} we get $\lambda\in]0,1]$ and $\ol\be_i=(1-\lambda)\,K_i$.
Letting $\ol\a_i:=\a_i+\lambda\,K_i$, we have $\ol\a_i=0=f_i(0)$.
\item If $\Delta\Theta_i<0$, there exist $\m\geq 0$ such that $\ol\be_i=\m\,K_i$, so that \eqref{GBlocalter} becomes $\Delta\Theta_i=-(\lambda+\m-1)\,K_i$, whence $\lambda+\m>1$.
If $\lambda\geq 1$, letting $\ol\a_i:=\a_i+K_i$, this time we have
$$ f_i(0)\geq\ol\a_i\,,\qquad  f_i(s_{i-1})=\Delta\Theta_i=\ol\a_i-\m\,K_i\leq\ol\a_i \,. $$
If $\lambda\in]0,1]$, letting $\ol\a_i:=\a_i+\lambda\,K_i$ we again have $\ol\a_i=0=f_i(0)$.
\item If $\Delta\Theta_i>0$, by \eqref{GBlocalter} we have $(1-\lambda)\,K_i>\ol\be_i$, hence $\lambda\in[0,1[$, so that we again let $\ol\a_i:=\a_i+\lambda\,K_i$.
\een
\par Finally, when $\a_i=0$, we have $\a_i=f_i(0)$, and we choose $\ol\a_i:=0$.
\par Recalling that $\ol s_i\in]s_{i-1},s_i]$, and setting $\Theta_{h,i}:=\ol\a_i+\Theta(s_{i-1})$, the proof of the Claim is completed as in the easier case $\gK\leq 0$ previously considered in Step~1.
\smallskip\par\noindent
{\sc Step~3:} We now check property~(b) in Proposition~\ref{PGB}. Denoting $\be_j:=\angle{\gt(s_{j})\gt_h(s_{j}-)}$, again by Lemma~\ref{LGB} and Theorem~\ref{TGB},
for each $j=1,\ldots,n$ we have
\beq\label{GBglobal} \Theta(s_j)-\Theta(s_{j-1})=\a_j-\be_j-\int_{U_j}\gi_{\GG_j}\,\gK\,dA\,. \eeq
In this formula, $\GG_j$ is the oriented closed curve given by the join of the arc of $\gc$ from $\gc(s_{j-1})$ to $\gc(s_{j})$ and the geodesic arc of $P_h$ from $\gc(s_{j})$ to $\gc(s_{j-1})$, and $\gi_{\GG_j}$ is the index of the curve $\GG_j$ on $\M$. Also, $U_j$ is
the open subset of $\M$ enclosed by the curve $\GG_j$.
\par Notice that by our construction, see Figure~\ref{GBpic}, we deduce that the index $\gi_{\GG_j}$ is well-defined and actually $\gi_{\GG_j}=\pm 1$ in the
interior of each component of $U_j$, whereas $\gi_{\GG_j}=0$ outside $U_j$.
Since moreover $\M$ is assumed smooth and compact, the Gauss curvature $\gK$ is uniformly bounded on $\M$. By \eqref{GBglobal}, we thus get:
\beq\label{estK} \Bigl\vert\int_{U_j}\gi_{\GG_j}\,\gK\,dA\Bigr\vert\leq \int_{U_j}|\gK|\,dA\leq \Vert\gK\Vert_\ii\cdot\meas(U_j)<\ii\,.  \eeq
\par Now, for $j=1,\ldots,n-1$, the oriented turning angle of the polygonal $P_h$ at the corner point $\gc(s_j)$ is equal to $\a_{j+1}-\be_j$, whereas
by \eqref{widT} we correspondingly get:
$$ \Theta_{h,j+1}-\Theta_{h,j}=(\a_{j+1}-\be_j)-\lambda_{j+1}\int_{\wid U_{j+1}}\gK\,dA+\lambda_{j}\int_{\wid U_{j}}\gK\,dA-\int_{U_j}\gi_{\GG_j}\,\gK\,dA\,,  $$
where $\lambda_j\in[-1,1]$, by our Claim, and $\wid U_j=\emp$, if $\a_j=0$.
By \eqref{estK} we can thus estimate:
\beq\label{jumpest}
 |\a_{j+1}-\be_j|\leq |\Theta_{h,j+1}-\Theta_{h,j}|+\Vert \gK\Vert_\ii\cdot\bigl(\meas(\wid U_{j+1})+\meas(\wid U_{j})+\meas(U_{j}) \bigr)\,.  \eeq
\par We now observe that the Jumps of the piecewise constant function $\Theta_h$ are the turning angles $(\a_{j+1}-\be_j)$,
whereas by \eqref{widT}, the corresponding Jumps of the modified angle function $\wid\Theta_h$ are equal to $(\Theta_{h,j+1}-\Theta_{h,j})$.
By summating on $j=1,\ldots,n-1$ in \eqref{jumpest}, and using that $\wid U_j\sb U_j$ for each $j$, we then infer:
$$ \Var(\Theta_h)\leq\Var(\wid\Theta_h)+\e_h\,,\qquad \e_h:=3\,\Vert \gK\Vert_\ii\cdot\sum_{j=1}^n\meas(U_j)\,. $$
\par Finally, by the uniform convergence of $P_h$ to $\gc$, we deduce that $\e_h\to 0$ as $h\to\ii$, whence property (b) in Proposition~\ref{PGB} holds true.
\epff
\section{Curves into Riemannian surfaces}\label{Sec:Riem}
In this section, we extend the previous results to the more general case of curves into Riemannian surfaces, i.e., 2-dimensional Riemannian manifolds $(\wid\M,g)$.
\par We assume that $\wid\M$ is smooth (at least of class $C^3$), closed, and compact.
Recall that we can always find a smooth isometric embedding $F:\wid\M\hookrightarrow\RN$ of $\wid\M$ into a surface $\M=F(\wid\M)$ immersed in
the $N$-dimensional Euclidean space, for some
$N\geq 4$.
Since the total intrinsic curvature of piecewise smooth curves involves the geodesic curvature
and the turning angles at corner points, we do not need $\wid\M$ to be oriented.
\adl\par\noindent
{\large\sc Total intrinsic curvature.} 
We first extend Definition~\ref{Dcurv}, by saying that the {\em total intrinsic curvature} of any curve $\g$ in $\wid\M$ is
$$ \TC_{\wid\M}(\g):=\lim_{\e\to 0^+}\sup\{ \gk_{\wid\M}(\wid P)\mid \wid P\in \Sigma_\e(\g)\} $$
where $ \Sigma_\e(\g)$ is the class of polygonals $\wid P$ in $\wid\M$ inscribed in $\g$ and with modulus $\m_\g(\wid P)<\e$, and $\gk_{\wid\M}(\wid P)$ is the rotation of $\wid P$,
both modulus and rotation being defined as in the case of surfaces $\M$ in $\RN$.
\adl\par\noindent
{\large\sc Results.} We extend the representation formula in Theorem~\ref{Tequality}, by the following:
\bt\label{TequalityS} Let $\wid\M$ be any smooth, closed, and compact Riemannian surface. For
every rectifiable curve $\g$ in $\wid\M$ with finite total intrinsic curvature, we have
$$ \TC_{\wid\M}(\g)=\F(\gt) $$
where the energy functional $\F(\gt)$ is defined by $\eqref{Ft}$ in correspondence to the tangent indicatrix $\gt=\dot\gc$ of $\gc=F\circ\g$, and $F$ is any isometric embedding of $\wid\M$ as above. \et
\par In order to prove Theorem~\ref{TequalityS}, we shall first introduce geodesic polar coordinates, and write the local expression \eqref{kgloc} of the geodesic curvature of a smooth curve $\g$ in $\wid\M$.
It turns out that length, angles and geodesics are preserved by isometries.
Letting then $\gc:=F\circ \g$, we shall compute the geodesic curvature $\kk_g$ of $\gc$ in $\M:=F(\wid\M)$, an immersed surface in $\RN$, showing that $\kk_g$
agrees with the intrinsic local expression \eqref{kgloc}, and hence that the latter does not depend on the choice of isometric embedding.
In a similar way, we will check that the rotation of a polygonal $\wid P$ in $\wid\M$ is an intrinsic notion.
\par As a consequence, we readily obtain:
\bp\label{PTCS} For any piecewise smooth curve $\g$ in $\wid\M$, we have
$$ \TC_{\wid\M}(\g)=\TC_\M(\gc)\qquad \If\quad \gc:=F\circ\g $$
independently of the chosen isometric embedding $F$.
\ep
\par Moreover, all the previous results obtained for curves $\gc$ in surfaces $\M$ of $\RN$ extend to curves $\g$ in a Riemannian surface $(\wid\M,g)$.
In fact, it suffices to
work with $\gc=F\circ \g$ for any isometric embedding $F$, and to use standard arguments based on local geodesic coordinates and partition of unity.
\par For this purpose, we shall focus in particular on the validity of the compactness theorem~\ref{Tcomp}.
In fact, by a quick inspection it turns out that the fundamental inequality \eqref{traspest} is the unique point of the previous theory where we used non-intrinsic quantities.
\par On account of Proposition~\ref{PTCS} and Theorem~\ref{Tequality}, we thus conclude with the validity of Theorem~\ref{TequalityS}.
\adl\par\noindent
{\large\sc Geodesic polar coordinates.} Following e.g. \cite[Sec.~4.12]{Du}, on small open domains $U$
of $\wid\M$ homeomorphic to a disk, we introduce geodesic polar coordinates $ds^2=dr^2+g(r,\f)\,d\f^2$, where $g$ is a non-negative smooth function on $U$.
We shall denote by $f_{,r},$ $f_{,\f}$, $f_{,rr}$, $f_{,r\f}$, and $f_{,\f\f}$ the partial first and second derivatives of a function $f(r,\f)$ on $U$.
The coefficient $g$ of the Riemannian metric satisfies
\beq\label{g} \lim_{r\to 0}g=0\,,\quad\lim_{r\to 0}(\sqrt {g})_{,r}=1\qquad\fa\,\f \eeq
compare \cite[Sec.~4.6]{doC}.
Also, in coordinates the non-trivial Christoffel coefficients of the Levi-Civita connection $\nabla_g$ of the Riemannian metric are
\beq\label{Gamma} \GG^1_{22}=-\frac 12\,g_{,r}\,,\quad \GG^2_{12}=\GG^2_{21}=\frac 1{2g}\,g_{,r}\,,\quad \GG^2_{22}=\frac 1{2g}\,g_{,\f}\,. \eeq
\par Let $\g:I\to\wid\M$ be a smooth and regular curve parameterized by arc-length. Assume that $\g(\wid I)\sb U$ for some open interval $\wid I\sb I$.
Also, we choose the pole of the coordinates not lying on the trace $\g(\wid I)$ of the curve. Therefore, there exists a positive real constant $c$ such that
$g(r,\f)\geq c>0$ for every $(r,\f)\in \g(\wid I)$.
\par In coordinates, we thus have
$\g(s)=(r(s),\f(s))$ for some smooth functions $r(s)$ and $\f(s)$ satisfying $\lan \dot\g(s),\dot\g(s)\ran_g=\dot r^2+g(r,\f)\,\dot\f^2=1$ for every $s\in \wid I$.
Therefore, the unit tangent vector and unit conormal are $$ \dot\g=(\dot r,\dot\f)\,,\qquad \dot\g^\perp:=(-g^{1/2}\dot\f,g^{-1/2}\dot r)\,. $$
\par The acceleration vector $\nabla_{\dot\g}\dot\g$ can be written in components as $(\nabla_{\dot\g}\dot\g)^k=\ddot\g^k+\GG^k_{ij}\dot\g^i\dot\g^j$, for $k=1,2$, so that in the previous local coordinates we get
\beq\label{acc}
 (\nabla_{\dot\g}\dot\g)^1=\ddot r-\frac 12\,g_{,r}\,\dot\f^2\,,\qquad (\nabla_{\dot\g}\dot\g)^2=\ddot\f+\frac 1{g}\,g_{,r}\,\dot r\,\dot\f+\frac 1{2g}\,g_{,\f}\,\dot\f^2\,.
\eeq
We have $\lan \nabla_{\dot\g}\dot\g,\dot\g\ran_g=0$, whence $ \nabla_{\dot\g}\dot\g=\kk_g\,\dot\g^\perp$, where
$\kk_g:=\lan \nabla_{\dot\g}\dot\g,\dot\g^\perp \ran_g$ is the geodesic curvature of $\g$, so that $|\kk_g|=|\nabla_{\dot\g}\dot\g|_g$. This yields to the local expression:
\beq\label{kgloc} \ba{rl} \kk_g= & \ds\sqrt g\,\Bigl[ -\dot \f\,(\nabla_{\dot\g}\dot\g)^1
+\dot r\, (\nabla_{\dot\g}\dot\g)^2\Bigr]  \\
= &
\ds
\sqrt g\,\Bigl[(\dot r\,\ddot\f-\dot\f\,\ddot r)+\frac 12\,\Bigl( g_{,r}\,\dot\f^3+2\,\frac{g_{,r}}g\,\dot r^2\,\dot\f+\frac{g_{,\f}}g\,\dot r\,\dot\f^2 \Bigr)
\Bigr] \,. \ea \eeq
%
\bex If e.g. $\wid\M=\M=\Sph$ and $g(r,\f)=\sin^2r$, with $r=\theta$ and $\f=\vf$, using that
$$ \GG^1_{22}=-\sin\theta\cos\theta\,,\quad \GG^2_{12}=\GG^2_{21}=\cot\theta\,,\quad \GG^2_{22}=0$$
we recover the formula \eqref{geodcurvS} for $\kk_g$. \eex
\br\label{RGB} We also recall that if $\om$ denotes the angle between $\dot\g$ and the fixed direction $(1,0)$, we find
$$\tan\om=\sqrt g\,\frac{\dot\f}{\dot r}\,,\qquad \dot\om=\kk_g-(\sqrt g)_{,r}\,\dot\f\,. $$
Therefore, if the curve $\g$ parameterizes the positively oriented boundary of the smooth domain $U$, by Stokes theorem, compare \cite[Sec.~4.12]{Du}, one has
\beq\label{Gauss}
 \oint_{\pa U} (\sqrt g)_{,r}\,\dot\f\,ds=-\int_U\gK\,dA\,,\qquad \gK= -\frac 1{\sqrt g}\,(\sqrt g)_{,rr} \eeq
where $\gK$ is the Gauss curvature of $(\M,g)$, yielding to the local formula of Gauss-Bonnet theorem:
\beq\label{GB}
 \int_U\gK\,dA= 2\pi- \oint_{\pa U} \kk_g\,ds\,. \eeq
\er
{\large\sc Embeddings.} Given an isometric embedding $F:\wid\M\hookrightarrow\M\sb\RN$, we let $\ol g$ and $\ol\nabla$ denote the (Gaussian) metric and
(Levi-Civita) connection induced by the Euclidean metric of $\RN$ on $\M$. The pull-back of $\ol g$ and of $\ol\nabla$ through $F$
agree with the metric $g$ and Levi-Civita connection $\nabla_g$ on $\M$, respectively. Therefore, in local coordinates as above, writing $F=F(r,\f):U\to \RN$, we have
\beq\label{dotF}  F_{,r}\bullet F_{,r}=1\,,\quad F_{,r}\bullet F_{,\f}=0\,,\quad F_{,\f}\bullet F_{,\f}=g\,. \eeq
By computing the partial second derivatives, 
we thus obtain the six formulas for the scalar products in $\RN$
\beq\label{ddotF}
\ba{lll}  F_{,r}\bullet F_{,rr}=0\,, & F_{,r}\bullet F_{,r\f}=0\,, & \ds F_{,r}\bullet F_{,\f\f}=-\frac 12\,g_{,r}\,, \\
 F_{,\f}\bullet F_{,rr}=0\,, & \ds F_{,\f}\bullet F_{,r\f}=\frac 12\,g_{,r}\,, & \ds F_{,\f}\bullet F_{,\f\f}=\frac 1{2}\,g_{,\f}\,. \ea \eeq
\par Letting $\gc(s):=F\circ \g(s)$, where $s\in \wid I$, the unit tangent vector and conormal corresponding to $\dot\g$ and $\dot\g^\perp$ take the expression
\beq\label{ut}
\gt=\dot r\,F_{,r}+\dot\f\,F_{,\f}\,,\qquad  \gu=-g^{1/2}\dot \f\,F_{,r}+g^{-1/2}\dot r\,F_{,\f}\,. \eeq
The curvature vector of the curve $\gc$ in $\RN$ then becomes
\beq\label{curvextr}
\gk=\dot\gt=\ddot r\,F_{,r}+\ddot\f\,F_{,\f} + \dot r^2\,F_{,rr}+2\,\dot r\,\dot\f\,F_{,r\f} +\dot \f^2\,F_{,\f\f}\,. \eeq
We compute the geodesic curvature of $\gc$ in $\M$ through the formula $\kk_g:=\dot\gt\bullet\gu$, obtaining by \eqref{dotF} and \eqref{ddotF}
$$ \ba{rl} \kk_g = & \ds -g^{1/2}\,\dot\f\,\Bigl(\ddot r+ \dot\f^2\,\Bigl(-\frac 12\,g_{,r} \Bigr) \Bigr)
 + g^{-1/2}\,\dot r\,\Bigl(g\,\ddot \f+ 2\,\dot r\,\dot\f\, \Bigl(\frac 12\,g_{,r} \Bigr)+ \dot\f^2\,\Bigl( \frac 1{2}\,g_{,\f}\Bigr)  \Bigr) \\
 = & \ds \sqrt g\,\Bigl[(\dot r\,\ddot\f-\dot\f\,\ddot r)+\frac 12\,\Bigl( g_{,r}\,\dot\f^3+2\,\frac{g_{,r}}g\,\dot r^2\,\dot\f+\frac{g_{,\f}}g\,\dot r\,\dot\f^2 \Bigr) \Bigr] \ea $$
which agrees with the local expression \eqref{kgloc} for the geodesic curvature of $\g$ in $\wid\M$.
\br\label{Rderiv} If $\g$ is a geodesic in $\wid\M$, the curve $\gc=F\circ\g$ is a geodesic in $\M$, whence the curvature vector $\dot \gt$ is orthogonal to both $F_{,r}$ and $F_{,\f}$. By \eqref{curvextr}, \eqref{dotF} and \eqref{ddotF} we have
$$ 0=\dot\gt\bullet F_{,r}=\ddot r-\frac12\,g_{,r}\,\dot\f^2\,, \qquad 0=\dot\gt\bullet F_{,\f}=g\,\ddot\f +g_{,r}\,\dot r\,\dot\f+\frac 1{2}\,g_{,\f}\,\dot\f^2 $$
and hence for a geodesic $\gc$ one recovers the local expressions of the equations $\nabla_{\dot\g}\dot\g=0$ from \eqref{acc}\,:
%
\beq\label{ddots}
 \ddot r=\frac12\,g_{,r}\,\dot\f^2\,,\qquad \ddot\f=-\frac 1{2g}\,\bigl( 2 g_{,r}\,\dot r\,\dot\f+{g_{,\f}}\,\dot\f^2\bigr)\,. \eeq
\er
{\large\sc Rotation of polygonals.} We now check that the rotation of a polygonal $\wid P$ in $\wid\M$ is an intrinsic notion.
Assume in fact that two geodesic arcs $\g_i$ of $\wid P$ meet at a point $(r_0,\f_0)$ in $U$.
Denoting by $(\dot r_i,\dot\f_i)$ the direction of the arc $\g_i$ at the point $(r_0,\f_0)$, where $i=1,2$, the rotation of $\wid P$ at $(r_0,\f_0)$ is equal to
$$\arccos\lan(\dot r_1,\dot\f_1),(\dot r_2,\dot\f_2)\ran_g=\arccos(\dot r_1\,\dot r_2+g(r_0,\f_0)\,\dot\f_1\,\dot\f_2)\,. $$
On the other hand, if $P=F(\wid P)$ is the corresponding polygonal in $\M=F(\wid\M)$,
the direction of the geodesic arc $\ol\g_i:=F\circ\g_i$ at the point $F(r_0,\f_0)$ is
$$\gv_i=\dot r_i\,F_{,r}(r_0,\f_0)+\dot\f_i\,F_{,\f}(r_0,\f_0) $$
and hence, using \eqref{dotF}, the corresponding rotation angle is
$$\arccos (\gv_1\bullet\gv_2)=\arccos(\dot r_1\,\dot r_2+g(r_0,\f_0)\,\dot\f_1\,\dot\f_2)\,. $$
\par Therefore, the rotation of $\wid P$ is equal to the rotation of $P$, i.e., $\gk_{\wid\M}(\wid P)=\gk_{\M}(P)$,
independently of the chosen isometric embedding $F:\wid\M\hookrightarrow\M\sb\RN$.
\adl\par\noindent
{\large\sc The compactness theorem.}
Going back to Theorem~\ref{Tcomp} on the $W^{1,1}$ compactness of the transport vector fields,
it turns out that the fundamental inequality \eqref{traspest} actually involves a constant factor $C_\M$ which depends on the surface $\M$, see Remark~\ref{Rest}.
Therefore, in the case of curves in a Riemannian surface $(\wid\M,g)$, the constant $C_\M$ definitely depends on the chosen embedding $F$.
\par However, since $\wid\M$ is assumed to be of class $C^3$ and compact, all the derivatives of $F$ up to the third order are equibounded on $U$, independently of the local chart on $\wid\M$.
Moreover, with the previous notation, we may and do assume that $g(r,\f)\geq c>0$ on $\g(\wid I)$, where the positive constant (that depends on the choice of the poles of the polar geodesic coordinates) is independent of the normal neighborhood of the partition of $\wid\M$, by the smoothness and compactness of $\wid\M$.
Therefore, if $\wid P$ is a polygonal of $\wid\M$ inscribed in $\g$, by choosing the modulus $\m_\g(\wid P)$ sufficiently small, it turns out that
$g(r,\f)\geq c>0$ for each $(r,\f)$ in $\wid P$.
%
%
Setting then $P:=F\circ \wid P$, the above properties imply that (outside the corner points):
\ben \item both the normal curvatures $\kk_n$ of the polygonals $P$ in $\M$, and their derivatives w.r.t. the arc-length parameter, are equibounded by a constant only depending on $\M:=F(\wid\M)$;
\item if $\gu$ is the unit conormal of $P$, parameterized by arc-length, then
both $|\dot\gu|$ and $|\ddot\gu|$ are equibounded by a constant only depending on $\M:=F(\wid\M)$, see Example~\ref{Edotu}. \een
\par
By the above construction, we deduce that our compactness result continues to hold.
\bex\label{Edotu}
We finally check that the local expressions of the arc-length derivatives $\dot\gu$ and $\ddot\gu$ of the unit conormal to the curve $\gc:=F\circ\g$ do not depend on the second order
derivatives of $r$ and $\f$, when $\g$ is a geodesic arc in $\wid\M$.
\par By formula \eqref{ut}, in fact, in general we obtain
$$ \dot\gu= a\,F_{,r}+ b\,F_{,\f}+c\,F_{,rr}+d\,F_{,r\f}+e\,F_{,\f\f} $$
where
$$  a:= -\frac 1{2\sqrt g}\,(g_{,r}\,\dot r+g_{,\f}\,\dot\f)\,\dot\f-\sqrt g\,\ddot\f\,,\qquad
 b:= -\frac 1{2 g^{3/2}}\,(g_{,r}\,\dot r+g_{,\f}\,\dot\f)\,\dot r+\frac 1{\sqrt g}\,\ddot r $$
and
$$ c:= -\sqrt g\, \dot r\,\dot\f\,,\qquad d:=\frac 1{\sqrt g}\,(\dot r^2-g\,\dot\f^2)\,,\qquad e:=\frac 1{\sqrt g}\,\dot r\,\dot\f\,. $$
When $\g$ is a geodesic in $\M$, using the formulas \eqref{ddots} we can rewrite the coefficients $a$ and $b$ as
$$ a=\frac 1{2\sqrt g}\,{g_{,r}}\,\dot r\,\dot\f\,,\qquad b=  \frac 1{2 g^{3/2}}\,\bigl( {g_{,r}}\,(g\,\dot\f^2-\dot r^2)-{g_{,\f}}\,\dot r\,\dot\f\bigr)\,. $$
\par Therefore, when computing the second derivative $\ddot\gu$, using again the formulas \eqref{ddots} it turns out that its local expression only depends on the first derivatives of $(r,\f)$ and on the partial derivatives of $g$ and $F$ up to the third order, where, we recall, $g(r,s)\geq c>0$ along the given geodesic arc $\g$, as required.
\eex
\section{Development of curves}\label{Sec:displ}
The original idea of parallel transport by Tullio Levi-Civita involves the concept of {\em development} of a curve on a surface. If e.g. $\M=\Sph$, it corresponds to drawing in a plane
the points of the trace of the oriented curve in $\Sph$ as the 2-sphere rolls without slipping or spinning in the plane, while staying tangent to the plane at the points of the curve.
The above construction implies that the scalar curvature of the developed curve on $\gR^2$ is equal to the modulus of the geodesic curvature of the given curve in $\Sph$, see Example~\ref{Edev}.
\par In this final section, we analyze the relationship between the definition of total intrinsic curvature and the notion of development of a smooth curve,
see Proposition~\ref{Pdev}.
We point out that similar arguments, based on considering iterations of the development of the ``complete tangent indicatrix", are proposed by Reshetnyak \cite{Re}
as a way to treat the ``curvatures" of an irregular curve in $\RN$.
\adl\par\noindent{\large\sc Development of curves.} Following e.g. \cite{doC}, if $\g:I\to\M$ is a regular, smooth, and simple curve on a surface $\M\sb\gR^3$, and $\dot\gn(s)\neq 0$, where, we recall, $\gn(s)$ is the unit normal $\gn(s):=\dot\g(s)/\Vert\dot\g(s)\Vert$, then the {\em envelope of the tangent planes} is the ruled surface $\Sigma$ parameterized by
$$X(s,v):=\g(s)+v\,\frac{\gn(s)\tim\dot\gn(s)}{|\dot\gn(s)|} $$
that in the case $\M=\Sph$ clearly becomes $X(s,v):=\g(s)+v\,\gu(s)$.
Around the trace of the curve, the ruled surface $\Sigma$ has zero Gauss curvature, and hence, by Minding's theorem, it is locally isometric to a planar domain.
The parallel transport of tangent fields $X(s)$ along the curve is the same, when considering $\g$ either as a curve on $\M$ or as a curve on $\Sigma$.
In particular, when $X(s)=\gt(s)$, one can use either local coordinates on $\M$ or on $\Sigma$ in order to obtain the geodesic curvature $\kk_g$ of the curve $\g$.
As a consequence, the parallel transport can be computed locally by pulling back the parallel transport along the development of the curve on the plane $\gR^2$, see \eqref{displ}.
\par Moreover, we can define a tubular neighborhood (a strip) $\Sigma$ of the envelope of the tangent planes
to $\M$ along $\g$, in such a way that $\Sigma$ is a surface with Gauss curvature equal to zero.
As a consequence, the total curvature $\TC_\Sigma(\g)$ of $\g$ as a curve in $\Sigma$ is well-defined, according to Definition \ref{Dcurv}, by taking
inscribed polygonals $\wid P$ in $\Sigma$ with modulus sufficiently small (according to the width of the strip $\Sigma$, which actually depends on the maximum of the modulus of the geodesic curvature of the curve).
\par By means of the same vertexes as for $\wid P$, we may correspondingly consider the polygonal $P$ in $\M$ inscribed in $\g$.
However, in general the rotation
of $P$ in $\M$ is different from the rotation of $\wid P$ in $\Sigma$, i.e.,
$$ \gk_\M(P)\neq \gk_\Sigma(\wid P)\,. $$
\par In fact, if e.g. $\g$ is a parallel of the 2-sphere $\M=\Sph$, and the vertexes of $P$ are taken at equidistant points along $\g$, then the angles between $\wid P$ and $\g$ are equal to the angles between the developed curve in $\gR^2$ and the corresponding polygonal, whence they are smaller than the angles between $P$ and $\g$, see Example~\ref{Edev}.
\adl\par\noindent{\large\sc A representation formula.} Notwithstanding, we shall see that the total curvature $\TC_\Sigma(\g)$ of $\g$ in the strip $\Sigma$ can
be computed by means of its development:
\bp\label{Pdev} Let $\g$ be a regular, smooth, and simple curve on a smooth surface $\M\sb\gR^3$, with $\dot\gn\neq 0$ everywhere.
We have:
$\TC_\Sigma(\g)=\ds\int_\g|\kk_g|\,ds\,. $ \ep
\par Now, for any smooth curve $\g$ as in Proposition~\ref{Pdev}, Theorem~\ref{Tdens} says that the total curvature $\TC_\M(\g)$ agrees with the integral on the right-hand side of the previous formula, whence we get:
$$ \TC_\M(\g)=\TC_\Sigma(\g)\,. $$
\par In particular, if $\{P_h\}\sb \calP_\M(\g)$ satisfies $\m_\g(P_h)\to 0$, and $\{\wid P_h\}$ is (for $h$ large enough) the corresponding sequence of inscribed polygonals in $\Sigma$, even if in general one has $\gk_\M(P_h)\neq \gk_\Sigma(\wid P_h)$, we conclude that
$$ \lim_{h\to\ii}\gk_\M(P_h)=\lim_{h\to\ii}\gk_\Sigma(\wid P_h)=\ds\int_\g|\kk_g|\,ds\,. $$
\bpff{\sc of Proposition \ref{Pdev}:} By a standard covering argument, we can reduce to the case in which the trace of $\g$ is contained in a normal neighborhood $U$ of $\Sigma$, and we equip $U$ with geodesic polar coordinates where the pole does not lay on the trace of the curve $\g$.
By the local formula \eqref{Gauss} for the Gauss curvature, using that $\gK=0$ on $\Sigma$ it turns out that the coefficient $g$ of the Riemannian metric on $\Sigma$ satisfies
$(\sqrt g)_{,rr}=0$. By using the limits \eqref{g}, this yields that $g(r,\f)=r^2$, compare \cite[Sec.~4.6]{doC}.
\par
We thus have $\g(s)=(r(s),\f(s))$ for some smooth functions $r(s)$ and $\f(s)$ satisfying $\dot r^2+r^2\,\dot\f^2=1$, where $r(s)\geq c>0$ for each $s\in I$.
As a consequence, the acceleration vector in \eqref{acc} takes the form
$$
 (\nabla_{\dot\g}\dot\g)^1=\ddot r-r\,\dot\f^2\,,\qquad (\nabla_{\dot\g}\dot\g)^2=\ddot\f+\frac 2{r}\,\dot r\,\dot\f $$
and the local expression \eqref{kgloc} of the geodesic curvature $\kk_g$ of $\g$ becomes:
$$\kk_g= r\,(\dot r\,\ddot\f-\dot\f\,\ddot r)+(1+\dot r^2)\,\dot\f\,. $$
\par By Remark~\ref{RGB}, the angle $\om$ between $\dot\g$ and the fixed direction $(1,0)$ in the vector bundle $T\Sigma$ satisfies
$$\tan\om= r\,\frac{\dot\f}{\dot r}\,,\qquad \dot\om=\kk_g-\dot\f\,, $$
whence $\ds\frac d{ds}(\om +\f)=\kk_g $.
By Proposition~\ref{PTheta}, this yields that $\Theta:=\om+\f$ agrees (up to an additive constant) with the angle of the parallel transport along $\g$.
Therefore, any curve $\wid\g$ in $\gR^2$ with unit tangent vector
$$
\gT(s)=(\cos\Theta(s),\,\sin\Theta(s))\,,\qquad \Theta(s):=\om(s)+\f(s)  $$
is such that its scalar curvature agrees with $|\kk_g|$. Moreover, using that $\cos(\arctan x)=(1+x^2)^{-1/2}$ and $\sin(\arctan x)=x\,(1+x^2)^{-1/2}$,
and that $(1+x^2)=\dot r^{-2}\,(\dot r^2+r^2\,\dot\f^2)=\dot r^{-2}$ when $x=r\,\dot\f/\dot r$, we infer:
$$
\ba{rl} \cos(\om+\f)= &  \dot r\,\cos\f-r\,\dot\f\,\sin\f=\frac d{ds}(r\,\cos\f) \\
\sin(\om+\f)= &  \dot r\,\sin\f+r\,\dot\f\,\cos\f=\frac d{ds}(r\,\sin\f)\,. \ea  $$ 
\par As a consequence, up to a rigid motion in $\gR^2$, the local expression of the developed curve $\wid\g$ is:
\beq\label{displ}
 \wid\g(s)=r(s)\cdot\bigl(\cos\f(s),\, \sin\f(s)\bigr)\,,\qquad s\in  I\,. \eeq
Denoting then by $\gT(s)=\dot{\wid\g}(s)$ the unit tangent vector to $\wid\g$, by the previous computation we have
$$|\dot\gT(s)|=|\dot\Theta(s)|=|\dot\om(s)+\dot\f(s)|=|\kk_g(s)|\qquad\fa\, s\in I $$
and hence we deduce that
\beq\label{TCwidg} \TC(\wid\g)=\int_I|\kk_g(s)|\,ds \eeq
where, we recall, $\kk_g$ is the geodesic curvature of the given curve $\g$ as a curve in $\M$.
\par Now, for any sequence $\{\wid P_h\}\sb\calP_\Sigma(\g)$, condition $\m_\g(\wid P_h)\to 0$ implies that $\gk_\Sigma(\wid P_h)\to\TC_\Sigma(\g)$.
Moreover, by the above computation it turns out that
if $\wih P_h$ is the polygonal in $\gR^2$ inscribed in $\wid\g$ and with vertexes corresponding to the vertexes of $\wid P_h$ in $\g$, then
$\gk_\Sigma(\wid P_h)=\TC(\wih P_h)$ for each $h$. Also, property $\m_\g(\wid P_h)\to 0$ yields that $\mesh(\wih P_h)\to 0$. Since $\wid\g$ is a planar curve, we infer that
$\TC(\wih P_h)\to\TC(\wid\g)$. In conclusion, we get
$$ \TC_\Sigma(\g)=\lim_{h\to\ii}\gk_\Sigma(\wid P_h)=\lim_{h\to\ii}\TC(\wih P_h)=\TC(\wid\g)$$
and Proposition~\ref{Pdev} holds true on account of formula \eqref{TCwidg}. \epff
\bex\label{Edev} Following Example~\ref{Eparal}, if $\M=\Sph$ and $\g=\gc_{\theta_0}$ is the parallel with constant co-latitude $\theta_0\in ]0,\p/2]$,
the geodesic polar coordinates on $\Sph$ give $g=\sin^2 r$, so that $r(s)\equiv\theta_0$ and $\f(s)=s/\sin\theta_0$, where $s\in[0,2\p\,\sin\theta_0]$. The geodesic polar coordinates on $\Sigma$ give instead $g=r^2$, whence $r(s)\equiv \tan\theta_0$ and $\f(s)=\cot\theta_0\cdot s$, where again $s\in[0,2\pi\,\sin\theta_0]$.
Therefore, according to \eqref{displ}, the corresponding developed curve $\wid\g$ in $\gR^2$ is the arc of a circle of radius $\tan\theta_0$ and length $2\pi\,\sin\theta_0$, i.e.,
$$ \wid\g(s)=\tan\theta_0\,\bigl(\cos(\cot\theta_0\cdot s), \, \sin(\cot\theta_0\cdot s)\bigr)\,,\qquad s\in [0,2\pi\,\sin\theta_0]\,. $$
The pointwise scalar curvature of $\wid\g$ is the reciprocal of the curvature radius of $\wid\g$, and hence it is equal to the pointwise geodesic curvature $\kk_g\equiv \cot\theta_0$ of the parallel $\gc=\gc_{\theta_0}$, whereas the total curvature of $\wid\g$ is equal to $2\pi\,\cos\theta_0$, i.e., to the total curvature $\TC_\Sph(\gc_{\theta_0})$ of the parallel.
\eex
\adl\par\noindent{\bf Acknowledgements.}
We thank the referee for his/her suggestions that helped us to improve the presentation.
The research of D.M. was partially supported by the GNAMPA of INDAM.
The research of A.S. was partially supported by the GNSAGA of INDAM.
\end{document}